\documentclass[dvips,opre,nonblindrev,copyedit]{InformsStorage}% DEFAULT, double-spaced

\usepackage{endnotes}
\let\footnote=\endnote

\usepackage{storagepaper-or,mathmacros-or}

\usepackage{natbib}
 \bibpunct[, ]{(}{)}{,}{a}{}{,}%
 \def\bibfont{\small}%
\TheoremsNumberedThrough     % Preferred (Theorem 1, Lemma 1, Theorem 2)
\ECRepeatTheorems

\EquationsNumberedThrough    % Default: (1), (2), ...
\MANUSCRIPTNO{} % When the article is logged in and DOI assigned to it,
\begin{document}
%%%%%%%%%%%%%%%%

%\setcounter{page}{1}
%\DOI{10.1287/xxxx.0000.0000}%

\RUNAUTHOR{Su and El Gamal}

\RUNTITLE{Limits on the Benefits of Energy Storage for Renewable Integration}

% \TITLE{Bundling Information Goods of Decreasing Value}
\TITLE{Limits on the Benefits of Energy Storage for Renewable Integration}

\ARTICLEAUTHORS{%
\AUTHOR{Han-I Su}
\AFF{Department of Electrical Engineering, Stanford University, Stanford, CA 94305, \EMAIL{hanisu@stanford.edu}} %, \URL{}}
\AUTHOR{Abbas El Gamal}
\AFF{Department of Electrical Engineering, Stanford University, Stanford, CA 94305, \EMAIL{abbas@ee.stanford.edu}}
} % end of the block

\ABSTRACT{%
The high variability of renewable energy resources presents significant challenges to the operation of the electric power grid. Conventional generators can be used to mitigate this variability but are costly to operate and produce carbon emissions. Energy storage provides a more environmentally friendly alternative, but is costly to deploy in large amounts. This paper studies the limits on the benefits of energy storage to renewable energy: How effective is storage at mitigating the adverse effects of renewable energy variability? How much storage is needed? What are the optimal control policies for operating storage?  
To provide answers to these questions, we first formulate the power flow in a single-bus power system with storage as an infinite horizon stochastic program. We find the optimal policies for arbitrary net renewable generation process when the cost function is the average conventional generation (environmental cost) and when it is the average loss of load probability (reliability cost). We obtain more refined results by considering the multi-timescale operation of the power system. We view the power flow in each timescale as the superposition of a predicted (deterministic) component and an prediction error (residual) component and formulate the residual power flow problem as an infinite horizon dynamic program. Assuming that the net generation prediction error is an IID process, we quantify the asymptotic benefits of storage. With the additional assumption of Laplace distributed prediction error, we obtain closed form expressions for the stationary distribution of storage and conventional generation. Finally, we propose a two-threshold policy that trades off conventional generation saving with loss of load probability. We illustrate our results and corroborate the IID and Laplace assumptions numerically using datasets from CAISO and NREL.
% Enter your abstract
}%

%\KEYWORDS{deterministic inventory theory; infinite linear programming duality;
%  existence of optimal policies; semi-Markov decision process; cyclic schedule}

 %\HISTORY{Received April 2012.}

\maketitle
%%%%%%%%%%%%%%%%%%%%%%%%%%%%%%%%%%%%%%%%%%%%%%%%%%%%%%%%%%%%%%%%%%%%%%

% Samples of sectioning (and labeling) in MNSC
% NOTE: (1) \section and \subsection do NOT end with a period
%       (2) \subsubsection and lower need end punctuation
%       (3) capitalization is as shown (title style).
%
%\section{Introduction.}\label{intro} %%1.
%\subsection{Duality and the Classical EOQ Problem.}\label{class-EOQ} %% 1.1.
%\subsection{Outline.}\label{outline1} %% 1.2.
%\subsubsection{Cyclic Schedules for the General Deterministic SMDP.}
%  \label{cyclic-schedules} %% 1.2.1
%\section{Problem Description.}\label{problemdescription} %% 2.

% Text of your paper here

\section{Introduction}

The rapid increase in the world demand for electricity
%%%(EIA 2011, Figure 72)
\cite[Figure~72]{EIA2011} 
coupled with the need to reduce the high carbon emissions due to electric power generation from fossil fuel
%%%(EPA 2011, Table 3-7)
\cite[Table~3-7]{EPA2011}
are driving a dramatic increase in renewable energy generation from sources such as wind, solar, and hydro. The power generated from wind and solar, however, is intermittent and uncertain, which presents significant challenges to power system operation as the penetration of these sources increases
%%%(NREL 2010)
\citep{NREL2010}.
In the long timescale (weeks to hours), this variability causes power imbalances: When renewable generation falls short of meeting the demand, more conventional generation from combined-cycle combustion and gas turbines is needed, which increases power system operation cost and offsets some of the environmental benefits of renewable energy
%%%(Hart and Jacobson 2011)
\citep{Hart2011}; 
when renewable generation exceeds demand, the excess power generated must be curtailed. In the short timescale (minutes to seconds), the variability of renewable generation can lead to large frequency and voltage variations and higher loss of load probability.
%spinning reserves, should we use this phrase?
% SPINNING RESERVES MIGHT BE TOO RESTRICTED

In addition to using conventional generation, renewable energy variability can be mitigated architecturally via geographical generation diversity
%%%(NREL 2010)
\citep{NREL2010} 
and renewable resource diversity
%%%(Li et al. 2009)
\citep{Li2009}, 
and operationally using demand-response
%%%(Kirby and Milligan 2010)
\citep{Kirby2010} 
and energy storage
%%%(Bitar et al. 2011, Denholm et al. 2010)
\citep{Denholm2010,Bitar2011}. 
In particular, energy storage can help in two quite different ways
%%%(EPRI 2010)
\citep{EPRI2010}. 
\begin{itemize}
\item In the long timescale, bulk energy storage systems, such as pumped hydroelectric storage and compressed air energy storage (CAES), can be charged by the excess renewable energy generation during off-peak hours and discharged during peak hours, hence reducing the need for additional conventional generation capacity, and renewable energy curtailment.
 
\item In the short timescale, fast-response energy storage systems, such as flywheels and batteries, can also help improve reliability beyond what fast-ramping generation can achieve because of their much faster response time.
\end{itemize}

This paper aims to establish the limits on the benefits of storage for mitigating renewable energy variability. How much can storage help reduce the need for conventional generation? How much can it help improve reliability? How much storage is needed to reap these benefits? What are the optimal control policies that achieve these limits? Satisfactory answers to these questions can help in architecting the smart grid as well as in operating it efficiently and reliably. 
Since we wish to establish limits on the benefits of storage rather than analyzing the operation of a particular power system with storage, we will ignore the fixed and operating costs of storage as well as explicit economic benefits such as arbitrage (e.g. see
%%%Eyer and Corey (2010)
\cite{Eyer2010}) 
throughout this paper. Under certain assumptions, we will find that most of the benefits can be achieved with only a modest amount of storage. 

Following is an outline of the rest of the paper:
\begin{itemize}

\item In Section~\ref{sec:single-bus}, we consider a single-bus power system with storage and a slotted-time model for load and energy generation. We formulate the power flow problem as an infinite horizon average-cost stochastic program in which the input is the net renewable generation (difference between the renewable power generated and the demand) and the controls are conventional generation and the storage charging and discharging operations. We consider two cost functions, the expected average conventional generation (environmental cost) and the expected average loss of load probability (reliability cost). We find the optimal policies for each of these cost functions for arbitrary net renewable generation process (Theorems~\ref{thm:policy-generation} and~\ref{thm:policy-lossofload}). The performance of these policies is illustrated using load data from 
%%%CAISO (2012)
\cite{CAISO} 
and simulated renewable wind power generation data from NREL
%%%(Potter et al. 2008)
\citep{Potter2008}. 

\item In Section~\ref{sec:residual}, we make progress toward quantifying the degree to which storage can help mitigate the impacts of renewable generation and the amount of storage needed. We consider the multi-timescale operation of the power system (day-ahead, hour-ahead, minutes-ahead, and real time). For each timescale, we view power flow in the single-bus power system with storage as the superposition of a predicted (deterministic) component and an error (residual) random component. We assume that the predicted component is balanced (with possibly a fixed offset) and formulate the residual power flow problem as an infinite horizon average-cost dynamic program with the net renewable generation prediction error as input and fast-ramping generation and storage as control variables. Assuming that the net renewable generation error is an IID process, we show that storage can reduce fast-ramping generation (relative to no storage) by a factor that approaches the storage round-trip inefficiency as its capacity becomes large (Proposition~\ref{prop:fastramp-asymp}). We further show that storage can reduce the average loss of load probability to zero as storage capacity becomes large (Proposition~\ref{prop:lossofload-asymp}). 

We then observe from the NREL simulated wind power dataset that in the short timescale (minutes), the wind power generation prediction error is close to Laplace distributed. Under this additional assumption, we obtain closed form expressions for the minimum average fast-ramping generation and the stationary distribution of the stored power sequence in some special cases (Propositions~\ref{prop:fastramp-lapl} and \ref{prop:fastramp-dist}). We show that most of the possible reduction in conventional generation can be achieved with relatively small storage capacity. We also show that the average loss of load probability can be reduced by an order of magnitude with small storage capacity.

\item The optimal policies we establish in Section~\ref{sec:single-bus} represent two extremes: The policy that minimizes the average conventional generation  always uses the stored energy ahead of conventional generation and never uses conventional generation to charge storage, while the policy that minimizes the average loss of load probability uses conventional generation ahead of stored energy and to keep storage as full as possible. In Section~\ref{sec:threshold}, we present a two-threshold policy that includes these two policies as special cases. We also show that this policy minimizes the one-period weighted sum of the aforementioned cost functions and is optimal for all the numerical examples we tried. Using this policy, we find a tradeoff between conventional generation capacity and storage capacity needed to achieve prescribed conventional generation consumption and loss of load probability.

\item In Section~\ref{sec:constrained} we generalize the policies established in Section~\ref{sec:single-bus} to the case where the storage charging and discharging rates are constrained (Theorems~\ref{thm:policy-generation-constr} and~\ref{thm:policy-lossofload-constr}). We find that most of the benefit to average conventional generation can be attained with relatively small charging and discharging rates, while higher such rates are needed for the average loss of load probability.
\end{itemize}

This paper is a significantly reorganized and expanded version of the conference paper
%%%(Su and El Gamal 2011)
\citep{2011}. 
There is a large body of previous work on energy storage. An overview of energy storage technologies and applications can be found in
%%%EPRI (2010), Roberts and Sandberg (2011)
\cite{EPRI2010, Roberts2011} 
and references therein. Prior related work to this paper include
%%%Chandy et al. (2010), Gayme and Topcu (2011), Oh (2011)
\cite{Chandy2010, Gayme2011, Oh2011}. 
Both
%%%Chandy et al. (2010)
\cite{Chandy2010} 
and
%%%Gayme and Topcu (2011)
\cite{Gayme2011}
assume that the net load is deterministic. 
%%%Chandy et al. (201)
\cite{Chandy2010} 
formulate a dynamic dc optimal power flow problem with energy storage as a convex program. The optimal policy is established explicitly for some special cases. 
%%%Gayme and Topcu (2011)
\cite{Gayme2011} formulate a dynamic ac optimal power flow problem with energy storage. The problem is shown to be non-convex in general, and sufficient conditions for strong duality are established.
%%% Oh (2011)
\cite{Oh2011} 
models renewable generation as a sequence of discrete random variables. An approximate stochastic programming method is proposed and illustrated via numerical examples. Concurrent and independent related work to this paper include
%%%Koutsopoulos et al. (2011)
\cite{Koutsopoulos2011} 
and
%%%ParandehGheibi et al. (2011)
\cite{ParandehGheibi2011} 
in which continuous-time models for the net renewable generation are considered with different problem formulations.

%--------------------
\section{Single-bus power system with storage}
\label{sec:single-bus}
% single-bus model

Consider the single-bus electric power system with storage depicted in Figure~\ref{fig:single-bus}, which consists of conventional generation, net renewable generation (difference between the renewable generation and the load), and energy storage. This power system may represent a transmission network with high renewable penetration, a distribution network with distributed renewable generators and energy storage devices, a microgrid not operated in the island mode where the power from the macrogrid acts as a fast-ramping generator
%%%(Lasseter et al. 2002)
\citep{Lasseter2002}, 
a wind farm with energy storage devices in which generation is acquired through an electricity market, or a stand-alone hybrid renewable energy system with battery storage
%%%(Bernal-Agustin and Dufo-Lopez 2009, Deshmukh and Deshmukh 2008))
\citep{Deshmukh2008,BernalAgustin2009}. The conventional generation may include base-load generators (coal-fired, hydro, and nuclear power plants), intermediate generators (combined-cycle combustion turbine), and peaking and fast-ramping generators (gas turbines).  The renewable generation may include wind and solar. The numerical results in this paper assume only wind power. The energy storage may include bulk energy storage (compressed air and pumped hydroelectric storage) and fast-response energy storage (flywheels and batteries).

\begin{figure}[h]
\FIGURE
{\tikzstyle{rectg}=[shape=rectangle,draw,minimum height=30pt,minimum width=70pt,text width=70pt, text centered]
\tikzstyle{rects}=[shape=rectangle,draw,minimum height=40pt,minimum width=70pt,text width=70pt, text centered]
\tikzstyle{rect}=[shape=rectangle,draw,minimum height=16pt,minimum width=40pt,text width=40pt, text centered]
\tikzstyle{multiline}=[text width=60pt, text centered]
\begin{tikzpicture}[font=\small]
\draw[very thick] (2,1.5) -- (2,4.6);
\node[style=rects,white,anchor=west] (sl) at (7.5,2.0) {Energy storage};
\node[style=rects,white,anchor=west] (sr) at (7.5,2.8) {Energy storage};
\node[style=rects,anchor=west] (s) at (7.5,2.4) {Energy storage $S_i \leq \Smax$};
\node[style=rect,anchor=west] (d) at (4,2.0) {discharge};
\node[style=rect,anchor=west] (c) at (4,2.8) {charge};
\node[style=rectg,anchor=east] (r) at (1,2.8) {Net renewable generation};
\node[style=rectg,anchor=west] (g) at (4,4) {Conventional generation};
\draw[-latex] (r) -- (2,2.8) node [midway,above] {$\D_i$};
\draw[-latex] (g) -- (2,4) node [midway,above] {$G_i \leq \Gmax$};
\draw[-latex] (d) -- (2,2.0) node [midway,above] {$D_i \leq \Dmax$};
\draw[-latex] (2,2.8) -- (c) node [midway,above] {$C_i \leq \Cmax$};
\draw[-latex] (c) -- (sr) node [midway,above] {$\ac C_i$};
\draw[-latex] (sl) -- (d) node [midway,above] {$D_i/\ad$};
\end{tikzpicture}}
{Single-bus power system with storage.
\label{fig:single-bus}}
{}
\end{figure}

% slotted-time
%--------------------
We assume a \emph{slotted-time} model for the dynamics of the power system, where time is divided into slots each of length $\tau$ hours and power is constant over each time slot. In the following, we introduce the needed definitions and the assumptions used throughout the paper.

\begin{itemize}
\item The  power supplied by the net renewable generation (difference between the renewable generation power and the load) in time slot $i =1,2,\ldots$ is denoted by $\D_i$. The sequence $\D_1,\D_2,\ldots$ is in general a random process.

\item We denote the total \emph{power capacity} of conventional generation by $\Gmax$ MW. The power supplied in time slot $i =1,2,\ldots$ is denoted by $G_i \leq \Gmax$.
\end{itemize}

We characterized energy storage by the following parameters:
\begin{itemize}
\item The \emph{energy storage capacity} $\tau \Smax$ MW-h is the maximum amount of energy that can be stored, where  $\Smax$ is referred to as the
\emph{power storage capacity}. Real-world energy storage devices cannot be completely discharged, and there is a limit on their minimum energy level. We use this minimum level as a reference and assume without loss of generality that it is equal to zero.
\item The \emph{stored power} at the beginning of time slot $i =1,2,\ldots$ is denoted by $S_i$ MW. 
\item The \emph{rated storage power conversion} $\Cmax$ MW is the maximum input (charging) power. The charging power at time $i =1,2,\ldots$ is denoted by $C_i \leq \Cmax$.
\item The \emph{rated storage output power} $\Dmax$ MW is the maximum output (discharging) power. The discharging power at time $i =1,2,\ldots$ is denoted by $D_i \leq \Dmax$.
\item 
The \emph{charging efficiency} $\ac \in (0,1)$ is the ratio of the charged power to the input power. The \emph{discharging efficiency} $\ad \in (0,1)$ is the ratio of the output power to the discharged power. The \emph{round-trip efficiency} therefore is $\a = \ac\ad$. 
\item The \emph{storage efficiency} is the fraction of retained power over a time slot. We assume throughout that the storage efficiency is very high compared to the round-trip efficiency and assume that it is equal to one. 
\end{itemize}

Using the above definitions, we can express the dynamics of the stored energy as
\begin{align*}
S_{i+1} & = S_i + \ac C_i - \frac{1}{\ad} D_i \quad \text{for}\ i =1,2,\ldots
\end{align*}
with the constraints $0 \leq S_i \leq \Smax$, $0 \leq C_i \leq \Cmax$, and $0 \leq D_i \leq \Dmax$, where $S_1 \in [0,\Smax]$ is given.

In Sections~\ref{sec:single-bus} through~\ref{sec:threshold}, we will assume unconstrained rated output power and power conversion of the energy storage, i.e.,
\begin{align}
\label{equ:unconstr}
\ac\Cmax &= \frac{1}{\ad}\Dmax = \Smax.
\end{align}
The general case is discussed in Section~\ref{sec:constrained}.

% balance constraints
%--------------------
We assume that the (negative) net renewable generation power is to be balanced as much as possible by conventional generation and stored power:
\begin{itemize}
\item If $\D_i \geq -\Gmax-\min\{\ad S_i,\Dmax\}$, then there is sufficient power capacity and the balance constraint must be satisfied, i.e., 
$G_i + D_i - C_i + \D_i \geq 0$. 
Note that if $G_i + D_i - C_i + \D_i > 0$, then there is excess generation. We assume that this excess generation is curtailed. 
\item If $\D_i < -\Gmax - \min\{\ad S_i,\Dmax\}$, then loss of load occurs. In this case, conventional generation is at its power capacity and storage is discharged at the rated storage output power, i.e.,  $G_i = \Gmax$, $C_i = 0$, and $D_i = \min\{\ad S_i,\Dmax\}$.
\end{itemize}

In the following two subsections, we formulate the power flow problem in the single-bus system as infinite horizon stochastic programs with two different cost functions. We establish the optimal control policies for both cost functions for arbitrary net renewable generation process. In Subsection~\ref{sec:single-bus-numerical}, we illustrate these policies using datasets from CAISO and NREL.

% single-bus problem i
%--------------------
\subsection{Minimizing average conventional generation}

The first cost function we consider is the expected long term average conventional generation. This is motivated by the need to reduce the carbon emissions of conventional generation. We seek to minimize this cost function by controlling the amount of conventional generation $G_i$, charging $C_i$, and discharging $D_i$ used in each time slot $i=1,2,\ldots$, where the triple $(G_i,C_i,D_i)$ is a function of the \emph{history} $H_i = (S_1,\D_1,S_2,\D_2,\ldots,S_{i-1},\D_{i-1},S_i)$. A control policy $\pi$ is a sequence of these triples, i.e., $\pi = \{(G_i,C_i,D_i):\, i = 1,2,\ldots\}$. A policy is said to be stationary if $\pi_i = \pi_j$ for all $i$ and $j$. We are now ready to define the first optimization problem.
\medskip

\noindent{Stochastic program I:} \emph{Minimizing average conventional generation}.
\begin{align*}
\text{minimize}\quad & \Jcg(\pi,S_1) = \limsup_{n \to \infty} \E\left[ \frac{1}{n} \sum_{i=1}^n G_i \right] \\
\text{subject to}\quad  
& S_{i+1} = S_i + \ac C_i - \frac{1}{\ad} D_i, \\
& 0 \leq S_i \leq \Smax, \qquad\ 0 \leq C_i, \\
& 0 \leq G_i \leq \Gmax, \qquad 0 \leq D_i, \\
& 0 \leq G_i - C_i + D_i+ \max\{\D_i, - \Gmax - \ad S_i \},
\end{align*}
where the expectation is over the net renewable generation sequence $\D_i$, $i=1,2,\ldots$. We denote the optimal policy by $\pg$.
\medskip

It turns out that a simple stationary policy is optimal for arbitrary net renewable generation process (including deterministic sequences).

\medskip
\begin{theorem}
\label{thm:policy-generation}
The optimal policy $\pg$ for stochastic program I is given in Table~\ref{tab:policy-generation}.
\end{theorem}
\medskip

The optimal policy is illustrated by the ``phase-diagram" in Figure~\ref{fig:policy-generation}. When the prediction error $\d \geq 0$, the optimal policy charges the storage using the excess renewable generation as much as possible. When the prediction error $\d < 0$, the storage is first discharged to compensate for as much of the renewable power deficit as possible. Conventional generation is then used to compensate for the remaining renewable generation deficit (if any). Thus, the optimal policy $\pg$ never charges the storage using conventional generation. This greedy policy is a consequence of the linearity of the cost function and the imperfect round-trip storage efficiency as will become clear in the proof.

\begin{table}[h]
\TABLE
{Optimal policy in Theorem~\ref{thm:policy-generation} and corresponding stored power dynamics.
\label{tab:policy-generation}}
{\begin{tabular*}{0.9\textwidth}{@{\extracolsep{\fill}}c|ccc|c}
$\pg_i$ & $\Gg_i$ & $\Cg_i$ & $\Dg_i$ & $S_{i+1}$ \\
\hline
$\displaystyle \frac{\Smax-S_i}{\ac} \leq \D_i$ & $0$ & $\displaystyle \frac{\Smax-S_i}{\ac}$ & $0$ & $\Smax$ \\
$\displaystyle 0 \leq \D_i < \frac{\Smax-S_i}{\ac}$ & $0$ & $\D_i$ & $0$ & $S_i + \ac \D_i$ \\
$-\ad S_i \leq \D_i < 0$ & $0$ & $0$ & $-\D_i$ & $S_i + \D_i /\ad$ \\
$- \Gmax - \ad S_i \leq \D_i < -\ad s_i$ & $-\D_i - \ad S_i$ & $0$ & $\ad S_i$ & $0$ \\
$\D_i < - \Gmax - \ad S_i$ & $\Gmax$ & $0$ & $\ad S_i$ & $0$ \\
\hline
\end{tabular*}}
{}
\end{table}

\begin{figure}[h]
\FIGURE
{\begin{tikzpicture}[scale=0.5,font=\small]
\node (s) at (9,0) {$S_i$};
\node (d) at (0,6) {$\D_i$};
\draw[-latex] (0,-5.5) -- (0,5.5);
\draw[-latex] (-0.5,0) -- (8.5,0);
\draw (8,5.5) -- (8,-5.5);
\draw (0,5) -- (8,0);
\node[anchor=east] at (-0.2,5) {$\Smax/\ac$};
\node at (6,3) {$\Smax$};
\node at (3,1.5) {$S_i + \ac \D_i$};
\node at (5,-0.7) {$S_i +\D_i/\ad$};
\node[anchor=east] at (-0.2,-2) {$-\Gmax$};
\draw (0,-2) -- (8,-5.2);
\node[anchor=west] at (8.2,-5.2) {$- \Gmax - \ad \Smax$};
\node at (3,-4.5) {$0$};
\draw (0,0) -- (8,-3.2);
\node[anchor=west] at (8.2,-3.6) {$-\ad \Smax$};
\node at (4,-2.5) {$0$};
\end{tikzpicture}}
{Illustration of the optimal policy $\pg$ in Theorem~\ref{thm:policy-generation}. The value in each region corresponds to the stored power at the end of slot $i$ when the stored power and prediction error at the beginning of this slot are $S_i$ and $\D_i$, respectively.
\label{fig:policy-generation}}
{}
\end{figure}

To prove the theorem, consider the finite horizon counterpart of stochastic program I. Define the cost-to-go function 
\begin{align*}
\Jg_k(\pi,H_k) = \E\left[ \left. \sum_{i=k}^n G_i \right| H_k \right].
\end{align*}
Let $\pi^*$ be the policy achieving $\inf_\pi \Jg_k(\pi,H_k)$.
For a fixed pair $(H_{k-1}, \D_{k-1})$, we will use the shorter notation $\Jg_k(\pi^*, S_k)$ in place of $\Jg_k(\pi^*, H_{k-1}, \D_{k-1}, S_k)$.
In the following lemma, we establish key properties of the cost-to-go function. The rest of the proof of Theorem~\ref{thm:policy-generation} is given in Appendix~\ref{app:single-bus}.

% lemma 1 for optimal policy i
%--------------------
\begin{lemma}
\label{lm:jg-upper}
The minimum cost-to-go function must satisfy the conditions:
\begin{enumerate}
\item $\Jg_k(\pi^*,S_k+c) \leq \Jg_k(\pi^*,S_k)$ for $0 \leq c \leq \Smax - S_k$, i.e., excess energy is not curtailed until storage is fully charged..

\item  $\Jg_k(\pi^*,S_k+c) \geq \Jg_k(\pi^*,S_k) - \ad c$ for $0 \leq c \leq \Smax - S_k$, i.e., the cost of using conventional generation to charge storage by $c$ MW is $c/\ac$, but the reduction in cost is at most $\ad c$. 
\end{enumerate}

\end{lemma}

\proof{}
To prove part 1, given the optimal policy $\pi^*$, we find another policy $\pi$ such that $\Jg_k(\pi,S_k+c) = \Jg_k(\pi^*,S_k)$. Let $H_i$ and $H_i^*$ be the history sequences under the policy $\pi$ and $\pi^*$, respectively. For fixed $(H_{k-1},\D_{k-1})$, consider $H_k = (H_{k-1},\D_k,S_k+c)$ and $H_k^* = (H_{k-1},\D_k,S_k)$. Let
\begin{align}
G_i & = G_i^*, \qquad
D_i = D_i^*, \qquad
C_i = \min\left\{ C_i^*, \frac{1}{\ac}\left( \Smax - S_i + \frac{1}{\ad} D_i^* \right) \right\}.
\label{equ:jg-upper-policy}
\end{align}
The policy $\pi$ is illustrated in Figure~\ref{fig:jg}\subref{fig:jg-upper}. By induction, the stored power under policy $\pi$ is always higher than that under policy $\pi^*$. Thus, by definition,
$\Jg_k(\pi^*,S_k) = \Jg_k(\pi,S_k + c) \geq \Jg_k(\pi^*,S_k + c).$

To prove part 2 of the lemma, given the optimal policy $\pi^*$, we find a policy $\pi$ such that $\Jg_k(\pi,S_k) \leq \Jg_k(\pi^*,S_k+c) + \ad c$. For fixed $(H_{k-1},\D_{k-1})$, consider $H_k = (H_{k-1},\D_k,S_k)$ and $H_k^* = (H_{k-1},\D_k,S_k+c)$. Let 
\begin{align}
C_i & = C_i^*, \qquad
D_i = \min\{ D_i^*, \ad(S_i + \ac C_i^*)\}, \qquad
G_i = \min\{ G_i^* + D_i^* - D_i, \Gmax \}.
\label{equ:jg-lower-policy}
\end{align}
The policy $\pi$ is illustrated in Figure~\ref{fig:jg}\subref{fig:jg-lower}. By induction, the stored power under policy $\pi$ is always lower than the stored power under policy $\pi^*$. The cost-to-go function of policy $\pi$ can be upper bounded as {\allowdisplaybreaks
\begin{align*}
\Jg_k(\pi,S_k) & \leq \E\left[ \left. \sum_{i=k}^n G_i^* + \ad( (S_i^*-S_i) - (S_{i+1}^* - S_{i+1}) ) \right| H_k \right] \\
& = \Jg_k(\pi^*,S_k+c) + \ad \E\left[ \left. (S_k^*-S_k) - (S_{n+1}^* - S_{n+1}) \right| H_k \right] \\
& \leq \Jg_k(\pi^*,S_k+c) + \ad c.
\end{align*}}
\endproof

\begin{figure}[h]
\FIGURE
{\subfloat[Part~1]{
\begin{tikzpicture}[scale=0.7,font=\footnotesize]
\draw[dashed,gray] (-2,-0.5) -- (8,-0.5);
\node at (-0.5,2) {time};
\node at (0.5,2) {$k$};
\node[anchor=east] at (-1.5,1) {$\pi^*$};
\node[anchor=east] at (-1.5,-2.1) {$\pi$};
\node[anchor=east] at (-0.2,1) {$S_k$};
\node[anchor=east] at (-0.2,-1.6) {$S_k+c$};
\node[anchor=east] at (-0.2,-2.6) {$c$};
\draw[pattern=north east lines] (0,0) rectangle (1,1);
\draw[thick] (0,1.6) -- (0,0) -- (1,0) -- (1,1.6);
\draw[fill] (0,-3) rectangle (1,-2.6);
\draw[pattern=north east lines] (0,-2.6) rectangle (1,-1.6);
\draw[thick] (0,-1.4) -- (0,-3) -- (1,-3) -- (1,-1.4);
\node at (2.5,2) {$k+1$};
\draw[pattern=north east lines] (2,0) rectangle (3,0.4);
\draw[thick] (2,1.6) -- (2,0) -- (3,0) -- (3,1.6);
\draw[fill] (2,-3) rectangle (3,-2.6);
\draw[pattern=north east lines] (2,-2.6) rectangle (3,-2.2);
\draw[thick] (2,-1.4) -- (2,-3) -- (3,-3) -- (3,-1.4);
\node at (3.75,0.8) {$\ldots$};
\node at (3.75,-2.2) {$\ldots$};
\node at (5,2) {$i$};
\draw[pattern=north east lines] (4.5,0) rectangle (5.5,1.4);
\draw[thick] (4.5,1.6) -- (4.5,0) -- (5.5,0) -- (5.5,1.6);
\draw[fill] (4.5,-3) rectangle (5.5,-2.8);
\draw[pattern=north east lines] (4.5,-2.8) rectangle (5.5,-1.4);
\draw[thick] (4.5,-1.4) -- (4.5,-3) -- (5.5,-3) -- (5.5,-1.4);
\node[anchor=west] at (5.7,-0.9) {Curtailing};
\draw[pattern=north east lines] (4.5,-1) rectangle (5.5,-0.8);
\node at (6.25,0.8) {$\ldots$};
\node at (6.25,-2.2) {$\ldots$};
%\draw[white] (4.5,2.5) rectangle (5.5,2.9);
\end{tikzpicture}
\label{fig:jg-upper}
} \qquad
\subfloat[Part~2]{
\begin{tikzpicture}[scale=0.7,font=\footnotesize]
\draw[dashed,gray] (-2,-1) -- (8.5,-1);
\node at (-0.5,2) {time};
\node at (0.5,2) {$k$};
\node[anchor=east] at (-1.5,0.9) {$\pi^*$};
\node[anchor=east] at (-1.5,-2) {$\pi$};
\node[anchor=east] at (-0.2,-2) {$S_k$};
\node[anchor=east] at (-0.2,1.4) {$S_k+c$};
\node[anchor=east] at (-0.2,0.4) {$c$};
\draw[pattern=north east lines] (0,-3) rectangle (1,-2);
\draw[thick] (0,-1.4) -- (0,-3) -- (1,-3) -- (1,-1.4);
\draw[fill] (0,0) rectangle (1,0.4);
\draw[pattern=north east lines] (0,0.4) rectangle (1,1.4);
\draw[thick] (0,1.6) -- (0,0) -- (1,0) -- (1,1.6);
\node at (2.5,2) {$k+1$};
\draw[pattern=north east lines] (2,-3) rectangle (3,-2.6);
\draw[thick] (2,-1.4) -- (2,-3) -- (3,-3) -- (3,-1.4);
\draw[fill] (2,0) rectangle (3,0.4);
\draw[pattern=north east lines] (2,0.4) rectangle (3,0.8);
\draw[thick] (2,1.6) -- (2,0) -- (3,0) -- (3,1.6);
\node at (3.75,0.8) {$\ldots$};
\node at (3.75,-2.2) {$\ldots$};
\node at (5,2) {$i$};
\draw[thick] (4.5,-1.4) -- (4.5,-3) -- (5.5,-3) -- (5.5,-1.4);
\draw[fill] (4.5,0) rectangle (5.5,0.2);
\draw[thick] (4.5,1.6) -- (4.5,0) -- (5.5,0) -- (5.5,1.6);
\node[anchor=east] at (4.3,-0.5) {$c'$};
\node[anchor=east] at (4.3,-3.46) {$\ad c'$};
\node[anchor=west] at (5.7,-0.5) {Discharging};
\node[anchor=west,text width=50pt] at (5.7,-3.3) {Conventional generation};
\draw[pattern=crosshatch dots] (4.5,-3.4) rectangle (5.5,-3.52);
\draw[fill] (4.5,-0.4) rectangle (5.5,-0.6);
\node at (6.25,0.8) {$\ldots$};
\node at (6.25,-2.2) {$\ldots$};
\end{tikzpicture}
\label{fig:jg-lower}
}
}
{Illustration for the proof of Lemma~\ref{lm:jg-upper}.
\label{fig:jg}}
{}
\end{figure}

% single-bus problem ii
%--------------------
\subsection{Minimizing average loss of load probability}
\label{sec:single-bus-lolp}
As we mentioned, energy storage can be used to improve power system reliability in the presence of renewable energy variation. As a measure of system reliability, we use the average loss of load probability. Define the loss of load cost at time $i$ as
\begin{align}
\label{equ:lossofload}
1_{\{\D_i < - \Gmax - \min\{\ad S_i, \Dmax\} \}} = 
\begin{cases} 
1 & \text{if } \D_i < - \Gmax - \min\{\ad S_i, \Dmax\} \\ 
0 & \text{otherwise.}
\end{cases}
\end{align}
We wish to find the control policy $\pi=\{(G_i,C_i,D_i):\, i=1,2,\ldots\}$ that minimizes the expected average loss of load cost.

\noindent{Stochastic program II:} \emph{Minimizing average loss of load}.
\begin{align*}
\text{minimize}\quad & \Jcl(\pi,S_1) = \limsup_{n \to \infty} \E\left[ \frac{1}{n} \sum_{i=1}^n 1_{\{\D_i < - \Gmax - \ad S_i \}} \right] \\
\text{subject to}\quad 
& S_{i+1} = S_i + \ac C_i - \frac{1}{\ad} D_i, \\
& 0 \leq S_i \leq \Smax, \qquad\ 0 \leq C_i, \\
& 0 \leq G_i \leq \Gmax, \qquad 0 \leq D_i, \\
& 0 \leq G_i - C_i + D_i+ \max\{\D_i, - \Gmax - \ad S_i \},
\end{align*}
where the expectation is over the net renewable generation sequence $\D_i$, $i=1,2,\ldots$. We denote the optimal policy by $\pl$.
\medskip

The optimal policy for this program is simply to keep the stored power as high as possible by using both excess renewable and conventional generation. The policy again holds for an arbitrary net renewable generation process.
%--------------------
\medskip
\begin{theorem}
\label{thm:policy-lossofload}
The optimal policy $\pl$ for stochastic program II is given in Table~\ref{tab:policy-lossofload}.  
\end{theorem}

\begin{table}[h]
\TABLE
{Optimal policy in Theorem~\ref{thm:policy-lossofload} and corresponding stored power.
\label{tab:policy-lossofload}}
{\begin{tabular*}{\textwidth}{@{\extracolsep{\fill}}c|ccc|c}
$\pl_i$ & $\Gl_i$ & $\Cl_i$ & $\Dl_i$ & $S_{i+1}$ \\
\hline
$\displaystyle \frac{\Smax-S_i}{\ac} \leq \D_i$ & $0$ & $\displaystyle \frac{\Smax-S_i}{\ac}$ & $0$ & $\Smax$ \\
$\displaystyle -\Gmax + \frac{\Smax-S_i}{\ac} \leq \D_i < \frac{\Smax-S_i}{\ac}$ & $\displaystyle \frac{\Smax-S_i}{\ac}-\D_i$ & $\displaystyle \frac{\Smax-S_i}{\ac}$ & $0$ & $\Smax$ \\
$\displaystyle -\Gmax \leq \D_i < -\Gmax + \frac{\Smax-S_i}{\ac}$ & $\Gmax$ & $\Gmax+\D_i$ & $0$ & $S_i+\ac(\Gmax+\D_i)$ \\
$\displaystyle -\Gmax - \ad S_i \leq \D_i < -\Gmax $ & $\Gmax$ & $0$ & $-\D_i-\Gmax$ & $\displaystyle S_i+ \frac{\Gmax+\D_i}{\ad}$ \\
$\D_i < -\Gmax - \ad S_i$ & $\Gmax$ & $0$ & $\ad S_i$ & $0$\\
\hline
\end{tabular*}}
{}
\end{table}

Since there is no conventional generation cost, in some cases the optimal policy charges the storage with conventional generation to minimize the loss of load probability as illustrated in the phase-diagram in Figure~\ref{fig:policy-lossofload},.  

\begin{figure}[h]
\FIGURE
{\begin{tikzpicture}[scale=0.5,font=\small]
\node (s) at (9.2,0) {$S_i$};
\node (d) at (0,6.2) {$\D_i$};
\draw[-latex,thick] (0,-5.5) -- (0,5.8);
\draw[-latex,thick] (-0.5,0) -- (8.8,0);
\draw (8,5.5) -- (8,-5.5);
\draw (0,5) -- (8,0);
\draw (0,3) -- (8,-2);
\node[anchor=east] at (-0.2,5) {$\Smax/\ac$};
\node at (6,3) {$\Smax$};
\node at (4,1.5) {$\Smax$};
\node[anchor=east] at (-0.2,3) {$-\Gmax + \Smax/\ac$};
\node at (7,-0.45) {$\Smax$};
\node[anchor=east] at (-0.2,0.8) {$S_i +\ac(\Gmax+\D_i)$}; 
\draw[-latex] (-0.3,0.8) -- (1,0.8);
\draw[-latex] (-0.3,0.8) -- (1.8,-0.8);
\node[anchor=west] at (8.5,-3) {$S_i + (\Gmax+\D_i)/\ad$};
\draw[-latex] (8.6,-3) -- (6.5,-3);
\draw (0,-2) -- (8,-2);
\node[anchor=east] at (-0.2,-2) {$-\Gmax$};
\draw (0,-2) -- (8,-5.2);
\node[anchor=west] at (8.2,-5.2) {$- \Gmax - \ad \Smax$};
\node at (3,-4.5) {$0$};
\end{tikzpicture}}
{Illustration of the optimal policy in Theorem~\ref{thm:policy-lossofload}.
\label{fig:policy-lossofload}}
{}
\end{figure}

% finite horizon, cost-to-go
To prove the theorem, we again consider its finite horizon counterpart and define the cost-to-go function as
\begin{align*}
\Jl_k(\pi,H_k) = \E\left[ \left. \sum_{i=k}^n 1_{\{\D_i < - \Gmax - \ad S_i \}} \right| H_k \right].
\end{align*}
Let $\pi^*$ be the policy achieving $\inf_\pi \Jl_k(\pi,H_k)$. 
For a fixed pair $(H_{k-1},\D_{k-1})$, we will used the shorter notation $\Jl_k(\pi,S_k)$ in place of $\Jl_k(\pi,H_{k-1},\D_{k-1},S_k)$.
The key step in the proof is to show
\begin{align}
\label{equ:jl-upper}
\Jl_k(\pi^*, S_k+c) \leq \Jl_k(\pi^*, S_k) \quad \text{for}\ 0 \leq c \leq \Smax - S_k.
\end{align}
This follows since the policy $\pi$ defined in~\eqref{equ:jg-upper-policy} always has higher stored power than the optimal policy $\pi^*$. Since the loss of load probability cannot increase if the amount of stored power is increased, we have
$\Jl_k(\pi^*, S_k) \geq \Jl_k(\pi,S_k + c) \geq \Jl_k(\pi^*,S_k + c).$
The rest of the proof of Theorem~\ref{thm:policy-lossofload} is given in Appendix~\ref{app:single-bus}.

% numerical results
%--------------------
\subsection{Numerical results}
\label{sec:single-bus-numerical}
\begin{figure}[h]
\FIGURE{
\footnotesize
\psfrag{xlabel}[c][c]{Time (hour)}
\psfrag{ylabel}[c][c]{Power (MW)}
\psfrag{legend-nrel}[r][r]{NREL wind dataset}
\psfrag{legend-caiso}[r][r]{CAISO load}
\includegraphics[width=0.45\textwidth,angle=270]{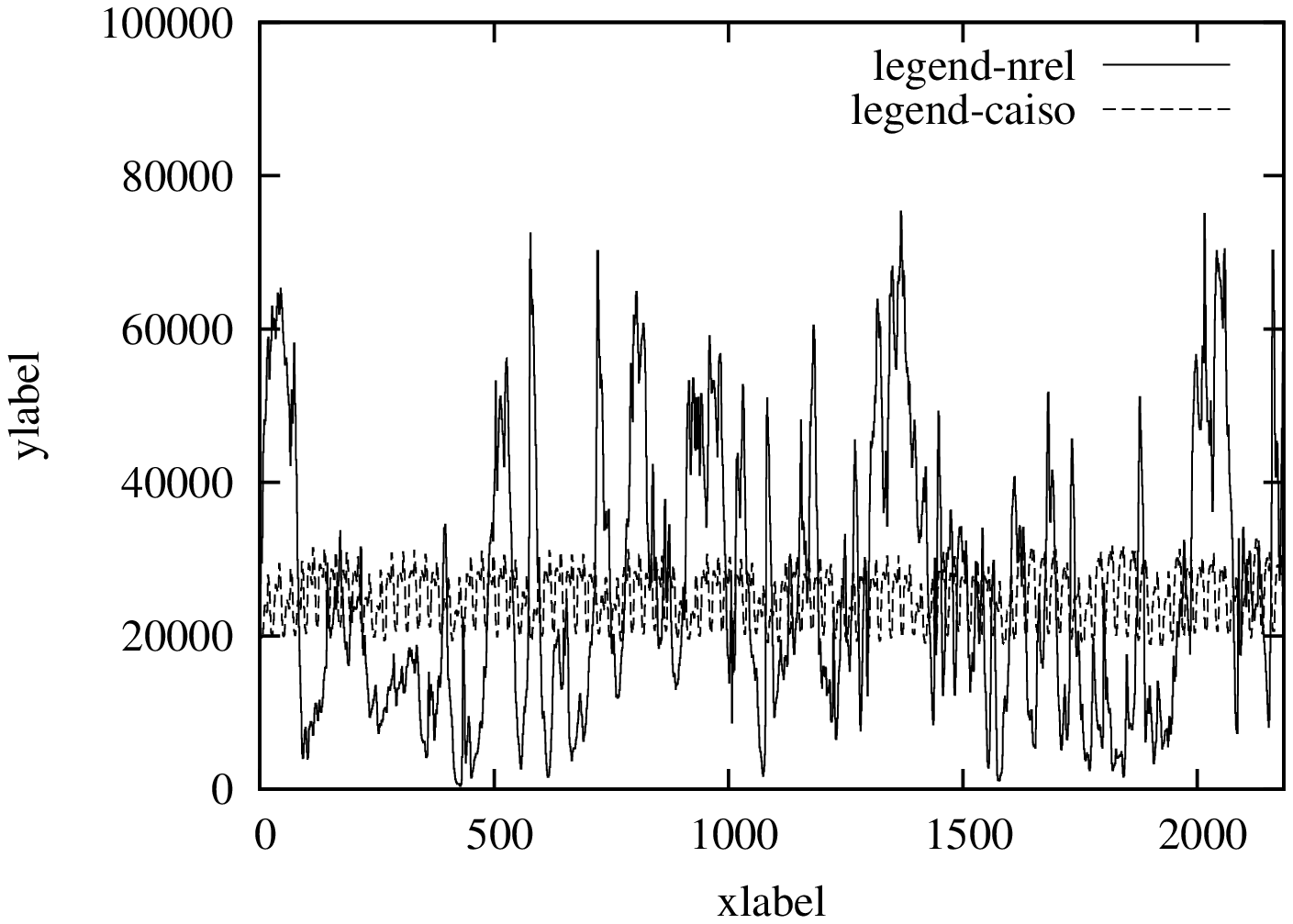}}
{The hourly average of the aggregate wind power of 3045 sites in California from NREL dataset and the actual load from CAISO dataset for three months in 2004. The average wind power is 25624 MW, and average load is 27232 MW.
\label{fig:sample-path}}
{}
\end{figure}

To investigate the reduction in conventional generation and loss of load probability attained by using energy storage, we use the simulated Western Wind Dataset from NREL
%%%Potter et al. 2008
\citep{Potter2008} 
and the actual load dataset from 
%%%%CAISO (2012)
\cite{CAISO}. The NREL dataset is based on numerical weather prediction (NWP) models. It attempts to recreate the potential wind power generation of more than 30\,000 sites with ten 3 MW turbines at each location in the western U.S from 2004 to 2006, and the wind power data is sampled every 10 minutes.
The CAISO load dataset includes the hourly load in California in 2004. 
% time series of wind power and load
%--------------------
Figure~\ref{fig:sample-path} plots the total simulated hourly power output of the wind turbines and the hourly load in California. 

\begin{figure}[h]
\FIGURE
{
\subfloat[]{
\footnotesize
\psfrag{xlabel}[c][c]{Power storage capacity (MW)}
\psfrag{ylabel}[c][c]{Conventional generation $\Jcg$ (MW)}
\psfrag{legend-nrel6}[r][r]{$\a = 60\%$}
\psfrag{legend-nrel8}[r][r]{$\a = 80\%$}
\includegraphics[width=0.35\textwidth,angle=270]{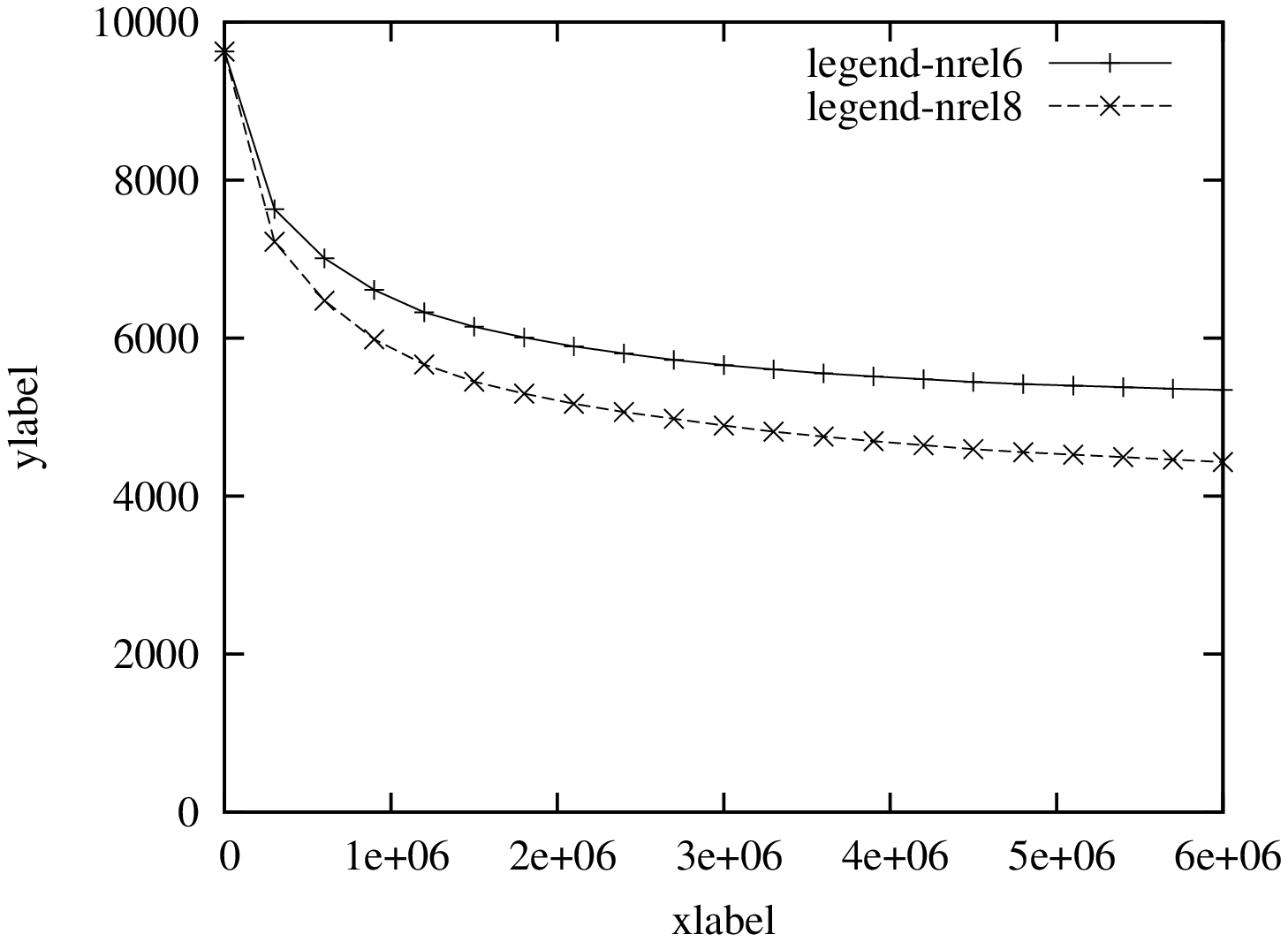}
\label{fig:generation-s}}
\subfloat[]{
\footnotesize
\psfrag{xlabel}[c][c]{Power storage capacity (MW)}
\psfrag{ylabel}[c][c]{Loss of load probability $\Jcl$}
\psfrag{legend-nrel6}[r][r]{$\a = 60\%$}
\psfrag{legend-nrel8}[r][r]{$\a = 80\%$}
\includegraphics[width=0.35\textwidth,angle=270]{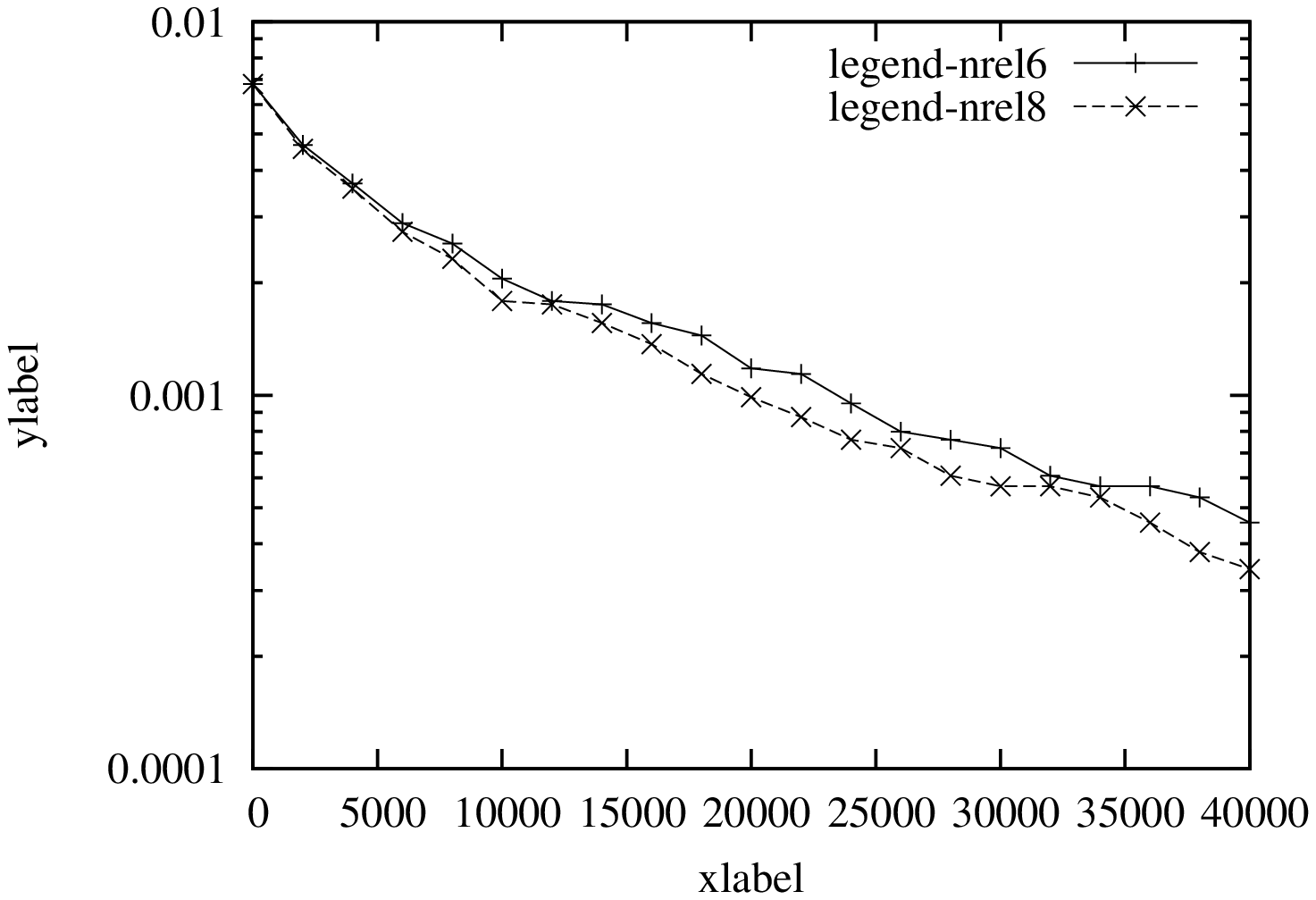}
\label{fig:lolp-s}}
}
{Figures~\subref{fig:generation-s} and \subref{fig:lolp-s} show the minimum average conventional generation and minimum loss of load probability, respectively, for NREL and CAISO data versus power storage capacity with round-trip efficiencies $\a = 60\%$ and $80\%$ and $\Gmax = 40\,000$ MW.
\label{fig:s}}
{}
\end{figure}

% conventional generation
%--------------------
Figure~\ref{fig:s}\subref{fig:generation-s} shows the average conventional generation under the optimal policy $\pg$ for several values of power storage capacities and round-trip storage efficiencies $\a = 60\%$ and $80\%$. 
As can be seen from the figure, conventional generation consumption can be reduced using storage by $44\%$ and $53\%$ for $\a = 60\%$ and $80\%$, respectively, when the storage capacity is 60 times the total rated power of the wind turbines, which is equivalent to  914 GW-h of storage. 
Note that $80\%$ of this reduction can be achieved by power storage capacity equal to 15 times the total rated power of the wind turbines, i.e., 228 GW-h.

% loss of load probability
%--------------------
Figure~\ref{fig:s}\subref{fig:lolp-s} plots the average loss of load probability under the optimal policy $\pl$ for several values of power storage capacities and round-trip storage efficiencies $\a = 60\%$ and $80\%$. 
For large $\Smax$, the reduction in the loss of load probability is roughly exponential in $\Smax$. We also see that
the average loss of load probability can be reduced by an order of magnitude with power storage capacity less than 2 standard deviations of the net renewable generation process.

\section{Residual power system}
\label{sec:residual}
% issues of single-bus formulation
%--------------------
In the previous section, we formulated the single-bus power system stochastic programs I and II and found the optimal policies $\pg$ and $\pl$ that minimize expected average conventional generation and average loss of load probability, respectively, for arbitrary net renewable generation process. These policies provide answers to the question concerning optimal storage control policies. Although answers to the other questions we seek to answer can be obtained numerically as illustrated in Subsection~\ref{sec:single-bus-numerical}, the results lump together vastly different types of conventional generation and storage resources that are deployed at very different timescales. Electric power systems are typically operated in multiple timescales as illustrated in Figure~\ref{fig:operation}:
\begin{itemize}
\item Day-ahead: Each day, an hourly prediction of the net load (the negative of the net renewable generation) for the next day is made. Base generation and bulk storage are scheduled to meet this prediction.

\item Hour-ahead: Each hour, a refined prediction of the net load in the next hour is made. Peaking generation and medium-response storage are scheduled to meet the difference (prediction error) between the day-ahead and hour-ahead prediction.

\item Minutes-ahead: Every few minutes, a prediction of the net load in the next few minutes is made. Fast-ramping generation and fast-response storage are scheduled to balance the difference between the minutes-ahead and hour-ahead prediction.

\item Real-time: The scheduled generation and storage are operated. The deviation of actual net load from the minutes-ahead prediction is matched by additional fast-ramping generation and fast-response storage.
\end{itemize}

\begin{figure}[h]
\FIGURE
{\begin{tikzpicture}[scale=0.9]
\draw[-latex] (0,0) -- (10,0);
\draw[-latex] (0,0) -- (0,5);
\node at (10.5,0) {time};
\node at (0,5.5) {MW};
\foreach \i in {2,3,...,8} {
\draw (\i,0.1) -- (\i,-0.1);
}
\node (h) at (5,-0.4) {Operating hour};
\foreach \i in {2,8} {
\draw[very thick] (\i,0.2) -- (\i,-0.2);
\draw[very thick] (\i,-0.3) -- (\i,-0.5);
\draw[-latex,thick] (h) -- (\i,-0.4);
}
\node[anchor=west] at (11.2,4.5) {\small Actual net load (real-time)};
\draw[thick] (10,4.5) -- (11,4.5);
\draw[smooth,thick] plot coordinates{
(1,0.8) (1.5,0.6) (2,0.9) (2.5,1.1) (3,1.6) (3.5,2.8) 
(4,2.4) (4.5,2.7) (5,3.4) (5.5,3.2) (6,4)   (6.5,3.5)
(7,3.8) (7.5,4.2) (8,4.3) (8.5,4.6) (9,4.5)};
\draw[latex-latex,thick] (6,3.5) -- (6,4);
\node[anchor=west] at (11.2,3) {\small Day-ahead};
\draw[dotted,very thick] (10,3) -- (11,3);
\draw[dotted,very thick] plot coordinates{
(1,1)   (3,2)   (7,2)   (9,3.5)};
\draw[latex-latex,dotted,very thick] (4.5,0) -- (4.5,2);
\node[anchor=west] at (11.2,3.5) {\small Hour-ahead};
\draw[dashed] (10,3.5) -- (11,3.5);
\draw[dashed] plot coordinates{
(1,0.5) (3,2.5) (7,2.5) (9,4.2)};
\draw[latex-latex,dashed,thick] (5.5,2) -- (5.5,2.5);
\node[anchor=west] at (11.2,4) {\small 10-minute-ahead};
\draw[densely dotted,very thick] (10,4) -- (11,4);
\draw[densely dotted,very thick] plot coordinates{
(1,0.9) (2,1.0) (3,1.5)
(4,2.6) (5,3.2) (6,3.5)
(7,3.7) (8,4.1) (9,4.4)};
\draw[latex-latex,densely dotted,very thick] (5,2.5) -- (5,3.2);
\end{tikzpicture}}
{Multi-timescale power grid operation.
\label{fig:operation}}
{}
\end{figure}

% residual model
%--------------------
Consider the multi-timescale power system operation. The day-ahead power flow prediction can be modeled in the same manner as the single-bus power system discussed in the previous section. In each subsequent timescale, the system can be decomposed into a scheduled (deterministic) part and a residual (random) part as depicted in Figure~\ref{fig:residual}. We assume that the scheduled generation and storage power balance the predicted net renewable generation power; hence we can model the \emph{residual power system} in exactly the same manner as the original system studied in Section~\ref{sec:single-bus}, except that the input is now the net renewable generation prediction error and the controls are fast-ramping generation and the charging and discharging of the fast-response energy storage. Unlike net renewable generation, which is in general a messy stochastic process, under good prediction, the prediction error process can be modeled as an IID zero-mean process $\D_i$, $i=1,2,\ldots$, with variance $\Var[\D_i] = \s^2$ (see Subsection~\ref{sec:residual-numerical} for numerical justification of this assumption). In Subsection~\ref{sec:residual-over-provision}, we consider the case of \emph{over provisioning} in which more power is scheduled than the predicted net renewable generation, that is, the mean of the prediction error $\E[\D_i]>0$.

\begin{figure}[h]
\FIGURE
{\tikzstyle{rectg}=[shape=rectangle,draw,minimum height=30pt,minimum width=70pt,text width=70pt, text centered]
\tikzstyle{rects}=[shape=rectangle,draw,minimum height=40pt,minimum width=70pt,text width=70pt, text centered]
\tikzstyle{rect}=[shape=rectangle,draw,minimum height=16pt,minimum width=40pt,text width=40pt, text centered]
\tikzstyle{multiline}=[text width=60pt, text centered]
\begin{tikzpicture}[font=\small]
\draw[very thick] (2,1.5) -- (2,4.6);
\node[style=rects,white,anchor=west] (sl) at (7.5,2.0) {Bulk energy storage};
\node[style=rects,white,anchor=west] (sr) at (7.5,2.8) {Bulk energy storage};
\node[style=rects,anchor=west,densely dashed] (s) at (7.5,2.4) {Bulk energy storage};
\node[style=rect,anchor=west,densely dashed] (d) at (4,2.0) {discharge};
\node[style=rect,anchor=west,densely dashed] (c) at (4,2.8) {charge};
\node[style=rectg,anchor=east,densely dashed] (r) at (1,2.8) {Net generation prediction};
\node[style=rectg,anchor=west,densely dashed] (g) at (4,4) {Scheduled generation};
\draw[-latex,densely dashed] (r) -- (2,2.8);
\draw[-latex,densely dashed] (g) -- (2,4);
\draw[-latex,densely dashed] (d) -- (2,2.0);
\draw[-latex,densely dashed] (2,2.8) -- (c);
\draw[-latex,densely dashed] (c) -- (sr);
\draw[-latex,densely dashed] (sl) -- (d);
\draw[very thick] (2,-2) -- (2,1.5);
\node[style=rects,white,anchor=west] (sl) at (7.5,-1.5) {Fast-response energy storage $S_i \leq \Smax$};
\node[style=rects,white,anchor=west] (sr) at (7.5,-0.7) {Fast-response energy storage $S_i \leq \Smax$};
\node[style=rects,anchor=west] (s) at (7.5,-1.1) {Fast-response energy storage $S_i \leq \Smax$};
\node[style=rect,anchor=west] (d) at (4,-1.5) {discharge};
\node[style=rect,anchor=west] (c) at (4,-0.7) {charge};
\node[style=rectg,anchor=east] (r) at (1,-0.7) {Net generation prediction error};
\node[style=rectg,anchor=west] (g) at (4,0.5) {Fast-ramping generation};
\draw[-latex] (r) -- (2,-0.7) node [midway,above] {$\D_i$};
\draw[-latex] (g) -- (2,0.5) node [midway,above] {$G_i \leq \Gmax$};
\draw[-latex] (d) -- (2,-1.5) node [midway,above] {$D_i \leq \Dmax$};
\draw[-latex] (2,-0.7) -- (c) node [midway,above] {$C_i \leq \Cmax$};
\draw[-latex] (c) -- (sr) node [midway,above] {$\ac C_i$};
\draw[-latex] (sl) -- (d) node [midway,above] {$D_i/\ad$};
\end{tikzpicture}}
{Decomposition of the total single-bus power system into a scheduled (predicted) power component and a residual power system component.
\label{fig:residual}}
{}
\end{figure}

% residual problem i
%--------------------
\subsection{Minimizing average fast-ramping generation}
\label{sec:residual-generation}
We first consider the role of fast-response energy storage in reducing the required fast-ramping generation. The formulation of this problem is the same as stochastic program I, but with the additional assumption that $\D_i$, $i = 1,2,\ldots$, is an IID process with zero-mean and variance $\Var[\D_i]=\s^2$. 
We refer to this new problem as \emph{dynamic program I}. 
Since the optimal policy $\pg$ in Theorem~\ref{thm:policy-generation} for stochastic program I holds for any net renewable generation process, it is optimal for dynamic program I. The IID assumption, however, allows us to provide some answers to the question of how much storage can help.

% asymptotic result
%--------------------

Consider the extreme case in which the fast-ramping generation capacity $\Gmax$ is unlimited. If there is no storage, i.e., $\Smax = 0$, then it is not difficult to see that the average cost is $\E[\D^-]$, where $x^- = \max\{-x, 0\}$. The minimum expected average fast-ramping generation for unlimited storage is given in the following. 

\medskip
\begin{proposition}
\label{prop:fastramp-asymp}
For unlimited $\Gmax$ and $\Smax$, the minimum average cost is $\Jcg(\pg,S_1) = (1-\a) \E[\D^-]$.
\end{proposition}
\medskip
The proof of this proposition is given in Appendix~\ref{app:residual}.

Comparing the average costs for no storage to unlimited power storage capacity, this proposition shows that storage can reduce the amount of needed fast-ramping generation (relative to no storage) in the limit by a factor equal to the round-trip storage inefficiency. This is not surprising because the IID zero-mean prediction error assumption implies that over the long term, the excess energy is roughly equal to the deficit. With infinite capacity and $\a = 1$, storage can compensate for almost all the variation in renewable generation. However, when $\a < 1$, it can compensate for at most this fraction of the variation and the rest needs to be compensated for by fast-ramping generation.

In exploring the wind generation prediction and prediction error data obtained using the NREL dataset (see Subsection~\ref{sec:residual-numerical} for details), we found that the first-order distribution of the 10-minute ahead prediction error $\D_i$ can be well-approximated by a Laplace$(\l)$ random variable with probability density function (pdf) $f_\D(\d) = e^{-\l|\d|}/2\l$. With this additional assumptions we can obtain the closed form expressions for the average fast-ramping generation and the stationary distributions of the stored power and fast-ramping generation in the following propositions.

\begin{proposition} 
\label{prop:fastramp-lapl}
The minimum expected average fast-ramping generation under the Laplace assumption is
\begin{align*}
\Jcg(\pg,S_1) & = \frac{1-e^{-\l\Gmax}}{2\l}\left( \frac{1-\a}{1-\a e^{-(1/\ac-\ad)\l\Smax/2}} \right).
\end{align*}
\end{proposition}
\medskip
The proof of this proposition is given in Appendix~\ref{app:residual}.

Using this result, we can obtain an answer to the question of how much fast-response storage is needed. Consider the derivative of the optimal cost function at storage capacity $\Smax$,
\begin{align*}
\left|\frac{\partial \Jcg(\pg,S_1)}{\partial \Smax} \right|
& = \frac{\ad(1-e^{-\l\Gmax})(1-\a)^2 e^{-(1/\ac-\ad)\l\Smax/2}}{4\left(1-\a e^{-(1/\ac-\ad)\l\Smax/2}\right)^2}.
\end{align*}
Since this derivative decreases close to exponentially, a small storage capacity is sufficient to achieve most of the reduction in the fast-ramping generation. For example, for a typical round-trip storage efficiency of 60\%--80\% (see
%%%Schainker (2004)
\cite{Schainker2004}), 
80\% of the possible reduction in the cost function can be achieved with power storage capacity less than 4 standard deviations of the prediction error.%, which is equivalent to 13.9 MW-h.

Note that under the optimal policy in Theorem~\ref{thm:policy-generation}, the stored power sequence is a homogeneous Markov process. In the following, we find the stationary distribution for this Markov process under the Laplace assumption. Furthermore, using this stationary distribution, we can find the distribution of the fast-ramping generation.

\begin{proposition}
\label{prop:fastramp-dist}
The cdf of the stationary distribution of the stored power under the optimal policy in Theorem~\ref{thm:policy-generation} and the Laplace assumption is
\begin{align*}
F_S(s) & = \frac{1 - 0.5(1+\a) e^{-(1/\ac-\ad)\l s/2}}{1-\a e^{-(1/\ac-\ad)\l\Smax/2}}
\end{align*}
for $0 \leq s < \Smax$, $F_S(s) = 0$ for $s < 0$, and $F_S(s) = 1$ for $s \geq \Smax$. The corresponding distribution of the fast-ramping generation is
\begin{align*}
F_G(g) & = 1 - \frac{1-\a}{2(1-\a e^{-(1/\ac-\ad)\l\Smax/2})} e^{-\l g}
\end{align*}
for $0 \leq g < \Gmax$, $F_G(g) = 0$ for $g < 0$, and $F_G(g) = 1$ for $g \geq \Gmax$.
\end{proposition}
\medskip
The proof of this proposition is given in Appendix~\ref{app:residual}.

Using the stationary distribution of the stored power sequence in Proposition~\ref{prop:fastramp-dist}, we can readily find the following expression for the expected average loss of load probability
\begin{align*}
\Jcl(\pg,S_1) &= \int_{-\infty}^{\Gmax} \frac{\l}{2} e^{\l\d} F_S\left(-\frac{\d+\Gmax}{\ad}\right) d\d % \\
% & = \frac{1}{1-\a e^{-(1/\ac-\ad)\l\Smax/2}} \left( \frac{1}{2}e^{-\l\Gmax} - \frac{1+\a}{2} \left(  \frac{\a}{1+\a} e^{-\l\Gmax}\right) \right) \\ &
= \frac{1}{2}e^{-\l\Gmax} \left( \frac{1-\a}{1-\a e^{-(1/\ac-\ad)\l\Smax/2}} \right).
\end{align*}

With no storage, the expected loss of load probability is $e^{-\l\Gmax}/2$. As $\Smax \to \infty$, $\Jcl(\pg,S_1) \to (1-\a) e^{-\l\Gmax}/2.$
Thus, under policy $\pg$, storage in the limit can reduce the expected average loss of load probability also by the round-trip storage inefficiency.

% residual problem ii
%--------------------
\subsection{Minimizing the loss of load probability}
\label{sec:residual-lolp}
We now consider the role of fast-response energy storage in reducing the average loss of load. The formulation of this problem is the same as stochastic program II with the additional assumption that $\D_i$, $i = 1,2,\ldots$, is an IID process with zero-mean and variance $\Var[\D_i]=\s^2$. 
We refer to this new problem as \emph{dynamic program II}. 
Since the optimal policy $\pl$ in Theorem~\ref{thm:policy-lossofload} for stochastic program II holds for any net renewable generation process, it is optimal for dynamic program II. 

%--------------------
In the following, we show that the benefit of storage to the loss of load probability is unbounded. 

\begin{proposition}
\label{prop:lossofload-asymp}
For unlimited storage capacity $\Smax$, if $F_\D$ and $\Gmax$ satisfy the conditions
\begin{align}
\limsup_{x\to\infty} x F_\D(-x) & \leq c, \label{equ:tail}\\
\E\left[ \a(\Gmax+\D)^+ - (\Gmax+\D)^- \right] & > 0, \label{equ:positive}
\end{align}
for some constant $c \geq 0$, then $\Jcl(\pl,S_1) = 0$.

\end{proposition}
The proof of this proposition is given in Appendix~\ref{app:residual}. 
Note that~\eqref{equ:tail} requires $f_\D$ with fast diminishing tail, and~\eqref{equ:positive} requires large enough $\Gmax$.

% Laplace remark
%--------------------
Unlike minimization of the expected average fast-ramping generation, we are not able to find closed form expressions for the optimal cost function or the stationary distributions under the Laplace assumption. However, it can be verified that for $\a > e^{-\l \Gmax}$, Laplace$(\l)$ satisfies the sufficient conditions in Proposition~\ref{prop:lossofload-asymp}, thus the expected average loss of load probability tends to zero as $\Smax \to \infty$.
Furthermore, in the following we show that the convergence rate of the loss of load probability is exponential in $\Smax$.

% Loss of load probability reducing rate with respect to Smax
%--------------------
\begin{proposition}
\label{prop:lossofload-rate}
For $\a > e^{-\l \Gmax}$, the exponent of the minimum expected average loss of load under the Laplace assumption decreases linearly as $\Smax$ increases, i.e.,
\begin{align*}
\gamma_{\mathrm{min}} \leq \lim_{\Smax \to \infty} \frac{\ln \Jcl(\pl,S_1)}{\Smax} \leq \gamma_{\mathrm{max}}
\end{align*}
for some constants $-\infty < \gamma_{\mathrm{min}} \leq \gamma_{\mathrm{max}} < 0$.
\end{proposition}
The proof is given in Appendix~\ref{app:residual}.

Note that when there is no fast-ramping generation, i.e., $\Gmax = 0$, the sufficient condition in Proposition~\ref{prop:lossofload-asymp} given by~\eqref{equ:positive} does not hold. However, the optimal policy reduces to a special case of the optimal policy in Theorem~\ref{thm:policy-generation}. Thus, storage can only reduce the expected loss of load probability by a factor no smaller than the round-trip storage inefficiency.

% Net generation prediction error with positive mean
%-------------------
\subsection{Over-provisioned net generation prediction error}
\label{sec:residual-over-provision}
Suppose that the net renewable generation prediction error $\D_i$ has mean $\E[\D_i] = \mu \geq 0$ and variance $\Var[\D_i] = \s^2$. In the following, we find the minimum expected average fast-ramping generation for unlimited $\Gmax$ and $\Smax$.

\medskip
\begin{proposition}
\label{prop:fastramp-over}
For unlimited $\Gmax$ and $\Smax$, the minimum average cost is
$\Jcg(\pg,S_1)  = \left( \E[\D^- - \a\D^+] \right)^+.$
\end{proposition}
\medskip

With no storage, i.e., $\Smax = 0$, the minimum average cost is $\Jcg(\pg,S_1) = \E[\D^-]$. Thus, storage can reduce the expected average fast-ramping generation $\Jcg$ by a factor of $1-\a\E[\D^+]/\E[\D^-]$ for small $\mu$ and reduce $\Jcg$ to $0$ for large $\mu$. 

For minimizing the expected loss of load probability, Proposition~\ref{prop:lossofload-asymp} still holds. Note that the minimum requirement on $\Gmax$ in~\eqref{equ:positive} decreases linearly in $\mu$.

% numerical results
%--------------------
\subsection{Numerical results}
\label{sec:residual-numerical}

\begin{figure}[h]
\FIGURE
{\footnotesize
\begin{minipage}{\textwidth}
\centering
\subfloat[Wind power]{
\psfrag{ylabel}[c][c]{Power (MW)}
\includegraphics[width=0.2\textwidth,angle=270]{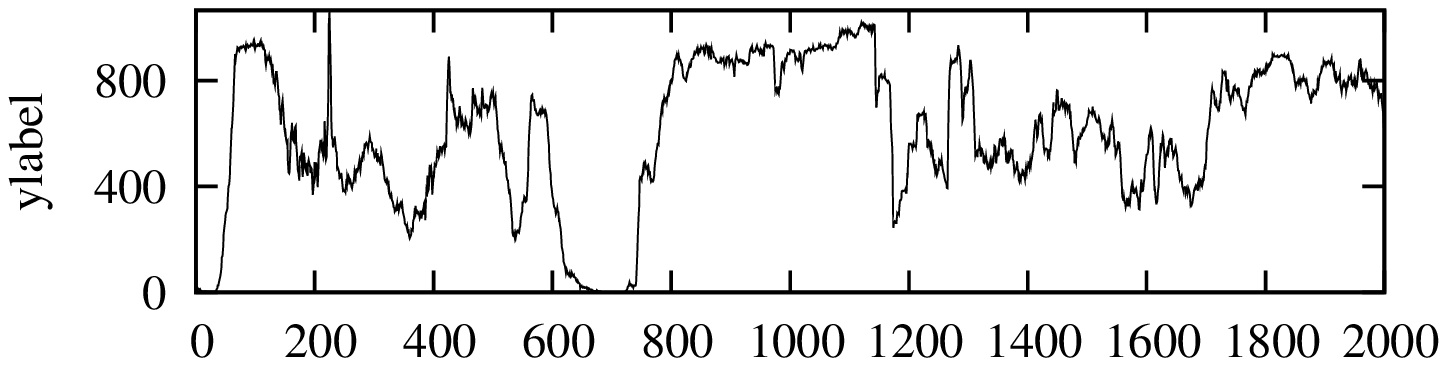}
\label{fig:wind-data}}\\
\subfloat[Prediction]{
\footnotesize
\psfrag{ylabel}[c][c]{Power (MW)}
\includegraphics[width=0.2\textwidth,angle=270]{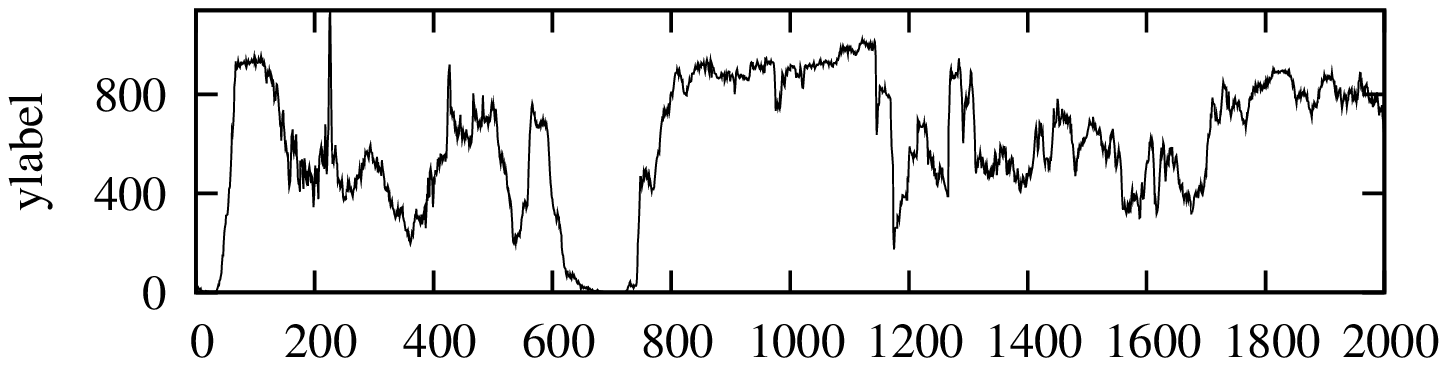}
\label{fig:wind-predict}}\\
\subfloat[Prediction error]{
\footnotesize
\psfrag{xlabel}[c][c]{Time index (10 minutes)}
\psfrag{ylabel}[c][c]{Power (MW)}
\includegraphics[width=0.2\textwidth,angle=270]{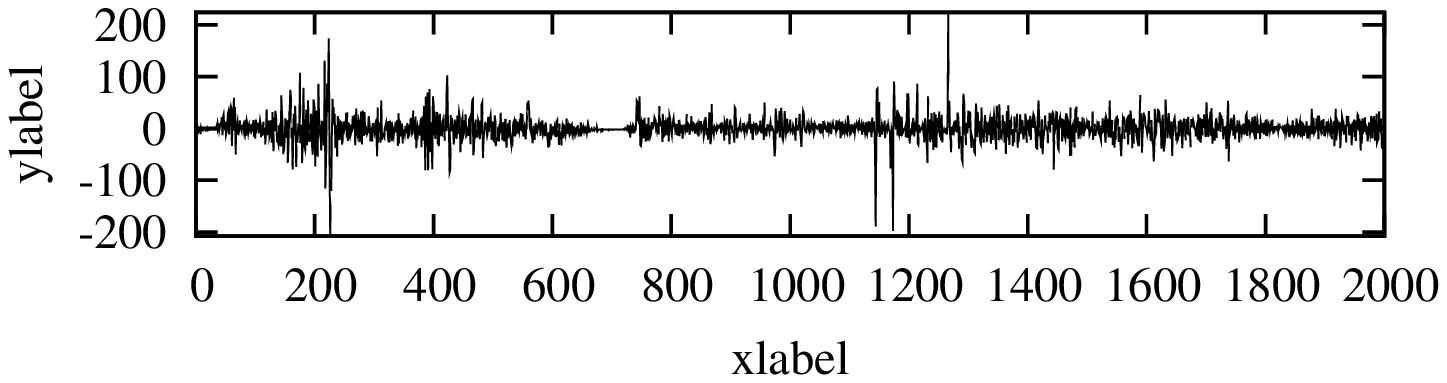}
\label{fig:wind-error}}
\end{minipage}
}
{Figure~\subref{fig:wind-data} shows the wind power over a two-week period. The average wind power is $560.26$ MW. Figure~\subref{fig:wind-predict} shows the 10-minute-ahead prediction given by the linear predictor based on the 6 samples in the past hour and optimized for the one-year data in 2004. The prediction error is shown in Figure~\subref{fig:wind-error}. The mean absolute value of the prediction error is $13.99$ MW, and the standard deviation of the prediction error is $20.88$ MW.
\label{fig:wind}}
{}
\end{figure}

First we provide numerical justifications for the IID and Laplace assumptions introduced in previous subsections. We use the NREL dataset for the 50 highest power density offshore wind sites in California. We assume that the variations in demand are much smaller than in the wind power, which is justified by the high penetration scenario assumed this paper (also see Figure~\ref{fig:sample-path}). Hence, we assume the net renewable generation error in our model is equal to the wind power prediction error. Since the dataset we use does not include forecast data, we use a simple linear predictor. 
Figure~\ref{fig:wind}\subref{fig:wind-data} plots the total power output of the 50 sites over a two-week period. 
Figures~\ref{fig:wind}\subref{fig:wind-predict} and~\ref{fig:wind}\subref{fig:wind-error} plot the 10-minute-ahead prediction and the prediction error sequences for the same two-week period, respectively.

\begin{figure}[h]
\FIGURE
{\footnotesize
\psfrag{xlabel}[c][c]{Power storage capacity $\Smax$ (MW)}
\psfrag{ylabel}[c][c]{Fast-ramping generation $\Jcg$ (MW)}
\psfrag{legend-nrel6}[r][r]{Wind dataset $\a = 60\%$}
\psfrag{legend-iid6}[r][r]{IID $\a = 60\%$}
\psfrag{legend-nrel8}[r][r]{Wind dataset $\a = 80\%$}
\psfrag{legend-iid8}[r][r]{IID $\a = 80\%$}
\includegraphics[width=0.45\textwidth,angle=270]{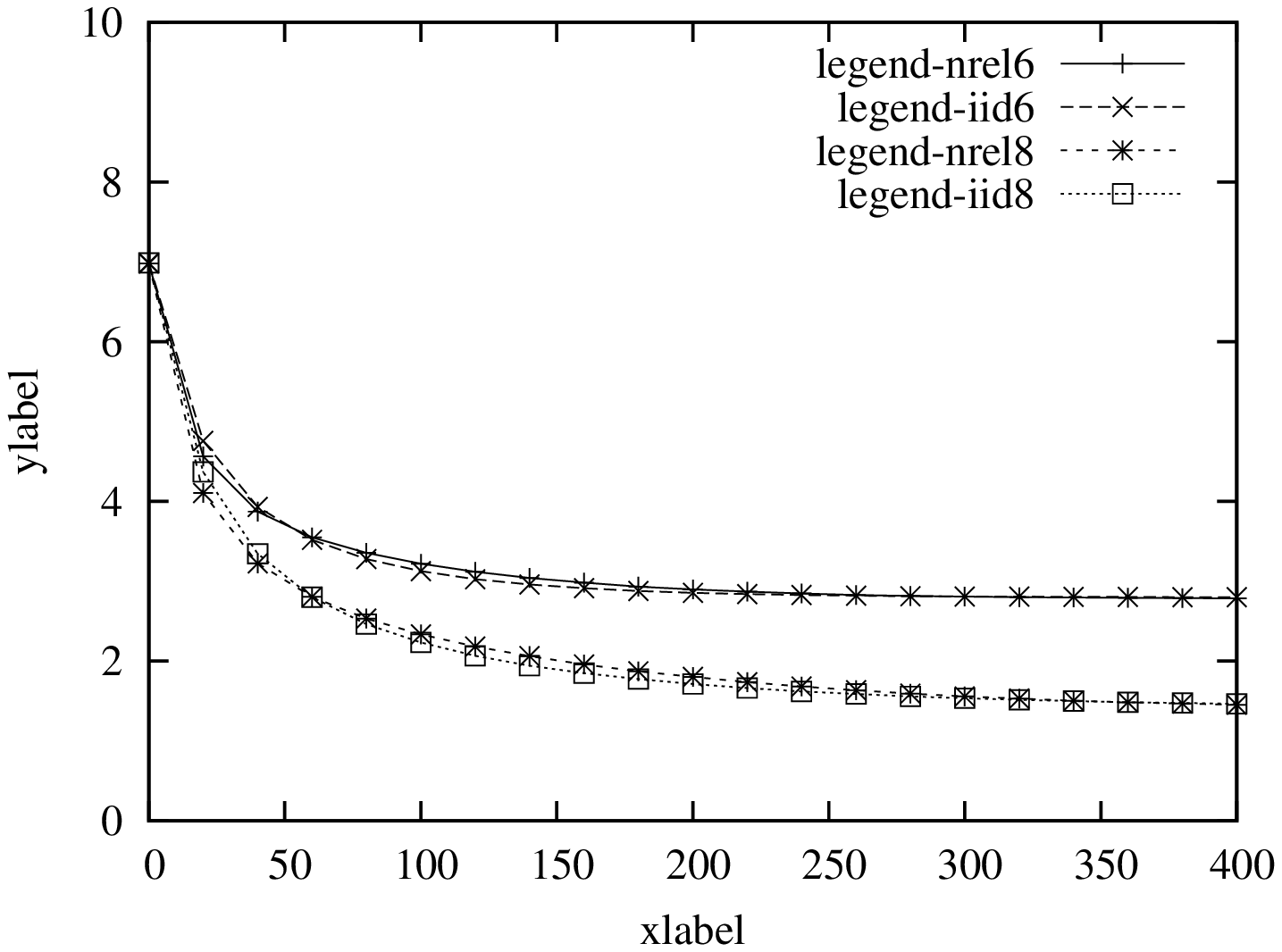}}
{The optimal expected average fast-ramping generation for the wind data versus power storage capacity for round-trip efficiencies $\a=60\%$ and $80\%$ and $\Gmax =160$ MW.
\label{fig:generation-iid-r}}
{}
\end{figure}

%As a sanity check of the IID assumption of the net renewable generation prediction error, Figure~\ref{fig:wind-cor} shows the autocovariance functions of the wind power and the 10-minute-ahead prediction error. We observe that the prediction error sequence is close to an uncorrelated sequence. ????????
To test the IID assumption, we generated a sequence of IID random variables distributed according to the empirical marginal distribution of the prediction error sequence from the NREL dataset. Figure~\ref{fig:generation-iid-r} compares the expected average fast-ramping generation costs for the NREL dataset and the IID sequence using policy $\pg$. The maximum absolute difference between costs for the NREL dataset and for the IID sequence normalized by the NREL dataset cost is less than 3\% and 6\% for $\a=60\%$ and $\a=80\%$, respectively.

\begin{figure}[h]
\FIGURE
{\footnotesize
\psfrag{xlabel}[c][c]{Prediction error (MW)}
\psfrag{ylabel}[c][c]{Probability density}
\psfrag{legendnrel}[r][r]{Prediction error}
\psfrag{legendlapl}[r][r]{Laplace($\l$) pdf}
\includegraphics[width=0.45\textwidth,angle=270]{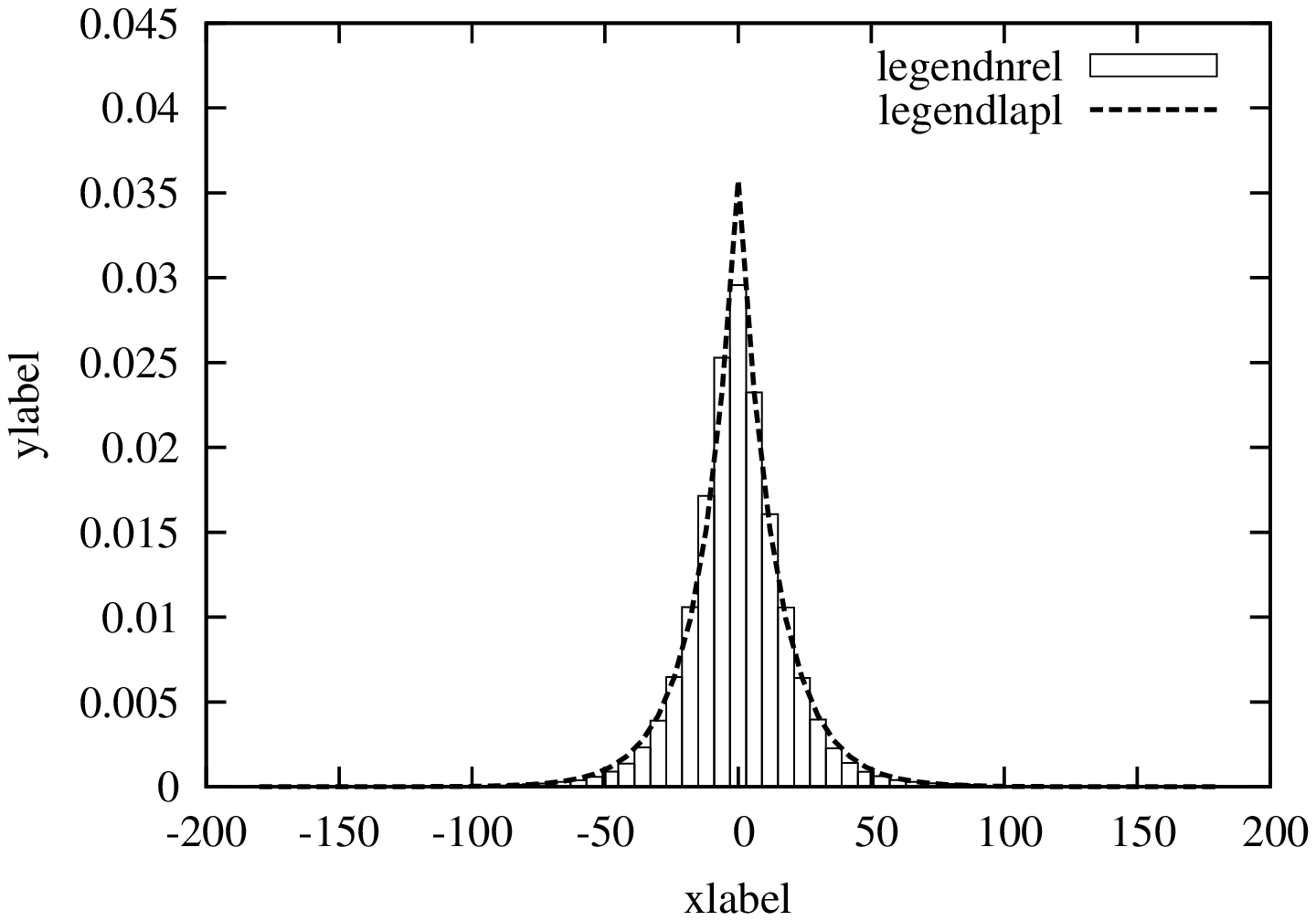}}
{The empirical pdf of the 10-minute-ahead wind power prediction error for three years versus the best fit Laplace($\l$) pdf with $1/\l = 13.99$.
\label{fig:wind-hist}}
{}
\end{figure}

To test the Laplace assumption, in Figure~\ref{fig:wind-hist} we compare the empirical pdf of the prediction error dataset to the best fit Laplace distribution. The maximum absolute difference between the empirical cdf of the prediction error and the Laplace cdf is $0.018$.
Hence, the assumption of IID Laplace distributed net renewable generation prediction error appears to be reasonable. In the following, we will further corroborate this assumption with the NREL data using the average cost results.

%\begin{figure}[h]
%\centering
%\footnotesize
%\psfrag{xlabel}[c][c]{Time index (10 minutes)}
%\psfrag{ylabel}[c][c]{Normalized autocovariance}
%\psfrag{legend-bulk}[r][r]{Wind power}
%\psfrag{legend-error}[r][r]{Prediction error}
%\includegraphics[width=0.52\textwidth,angle=270]{cor}
%\caption{The normalized autocovariance functions of the wind power and the 10-minute-ahead prediction error.}
%\label{fig:wind-cor}
%\end{figure}

% average generation v.s. Smax

\begin{figure}[h]
\FIGURE
{
\subfloat[]{
\footnotesize
\psfrag{xlabel}[c][c]{Power storage capacity $\Smax$ (MW)}
\psfrag{ylabel}[c][c]{Fast-ramping generation $\Jcg$ (MW)}
\psfrag{legend-nrel6}[r][r]{Wind dataset $\a = 60\%$}
\psfrag{legend-lapl6}[r][r]{Laplace($\l$) $\a = 60\%$}
\psfrag{legend-nrel8}[r][r]{Wind dataset $\a = 80\%$}
\psfrag{legend-lapl8}[r][r]{Laplace($\l$) $\a = 80\%$}
\includegraphics[width=0.37\textwidth,angle=270]{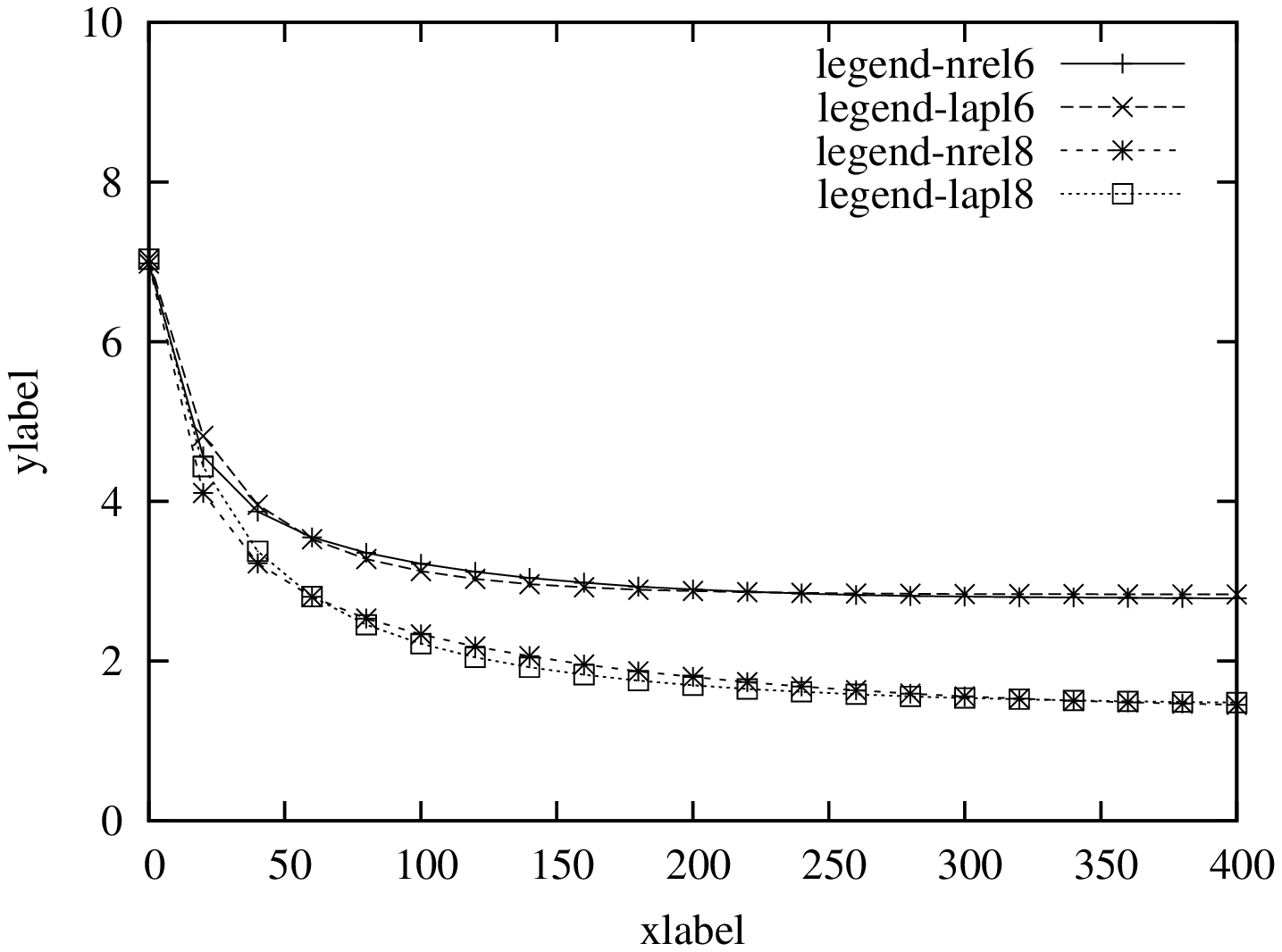}
\label{fig:generation-r}}
\subfloat[]{
\footnotesize
\psfrag{xlabel}[c][c]{Power storage capacity $\Smax$ (MW)}
\psfrag{ylabel}[c][c]{Loss of load probability $\Jcl$}
\psfrag{legend-nrel6}[r][r]{Wind dataset $\a = 60\%$}
\psfrag{legend-lapl6}[r][r]{Laplace($\l$) $\a = 60\%$}
\psfrag{legend-nrel8}[r][r]{Wind dataset $\a = 80\%$}
\psfrag{legend-lapl8}[r][r]{Laplace($\l$) $\a = 80\%$}
\includegraphics[width=0.37\textwidth,angle=270]{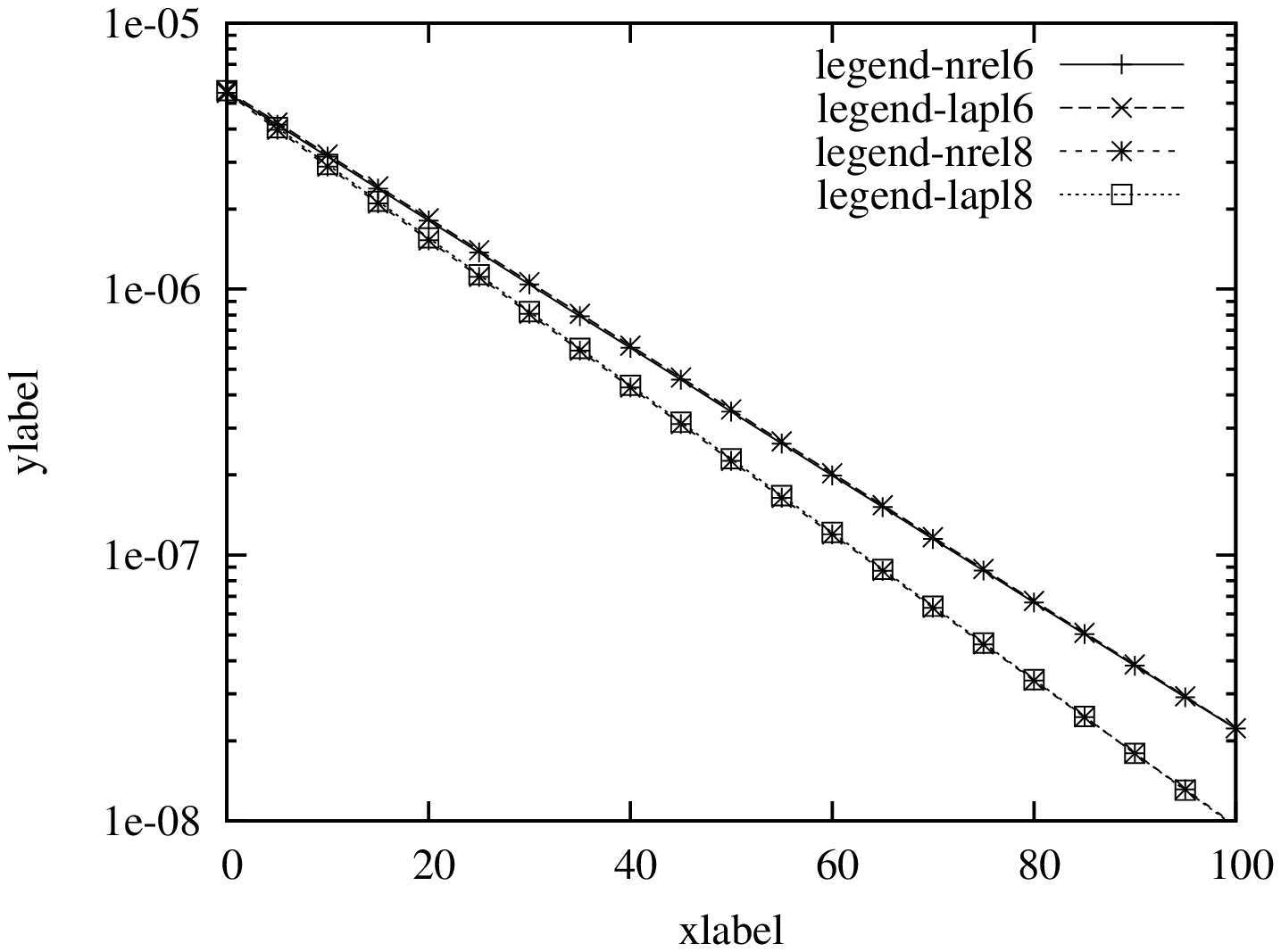}
\label{fig:lolp-r}}
}
{Figures~\subref{fig:generation-r} and \subref{fig:lolp-r} show the optimal expected average fast-ramping generation and the optimal loss of load probability, respectively, for the wind data versus power storage capacity for round-trip efficiencies $\a=60\%$ and $80\%$ and $\Gmax =160$ MW.
\label{fig:r}}
{}
\end{figure}

Now we illustrate the results in Subsections~\ref{sec:residual-generation} and \ref{sec:residual-lolp}. 
Figure~\ref{fig:r}\subref{fig:generation-r} compares the minimum average costs in Proposition~\ref{prop:fastramp-lapl} and using the three-year simulated wind data for various values of power storage capacities and round-trip storage efficiencies $\a=60\%$ and $80\%$.
The maximum absolute difference  between the theoretical and the simulated costs normalized by the theoretical cost is less than 6\% and 8\% for $\a=60\%$ and $\a=80\%$, respectively. Thus, the Laplace distribution appears to be an acceptable approximation of the simulated wind generation data from NREL.
Note that $80\%$ of the reduction for unlimited power storage capacity can be achieved with power storage capacity less than 4 standard deviations of the prediction error, which is equivalent to 13.9 MW-h.

\begin{figure}[h]
\FIGURE
{
\subfloat[Fast-response storage]{
\footnotesize
\psfrag{xlabel}[c][c]{Stored power (MW)}
\psfrag{ylabel}[c][c]{Probability density}
\psfrag{legendnrel}[r][r]{Wind dataset}
\psfrag{legendlapl}[r][r]{Laplace($\l$)}
\includegraphics[width=0.37\textwidth,angle=270]{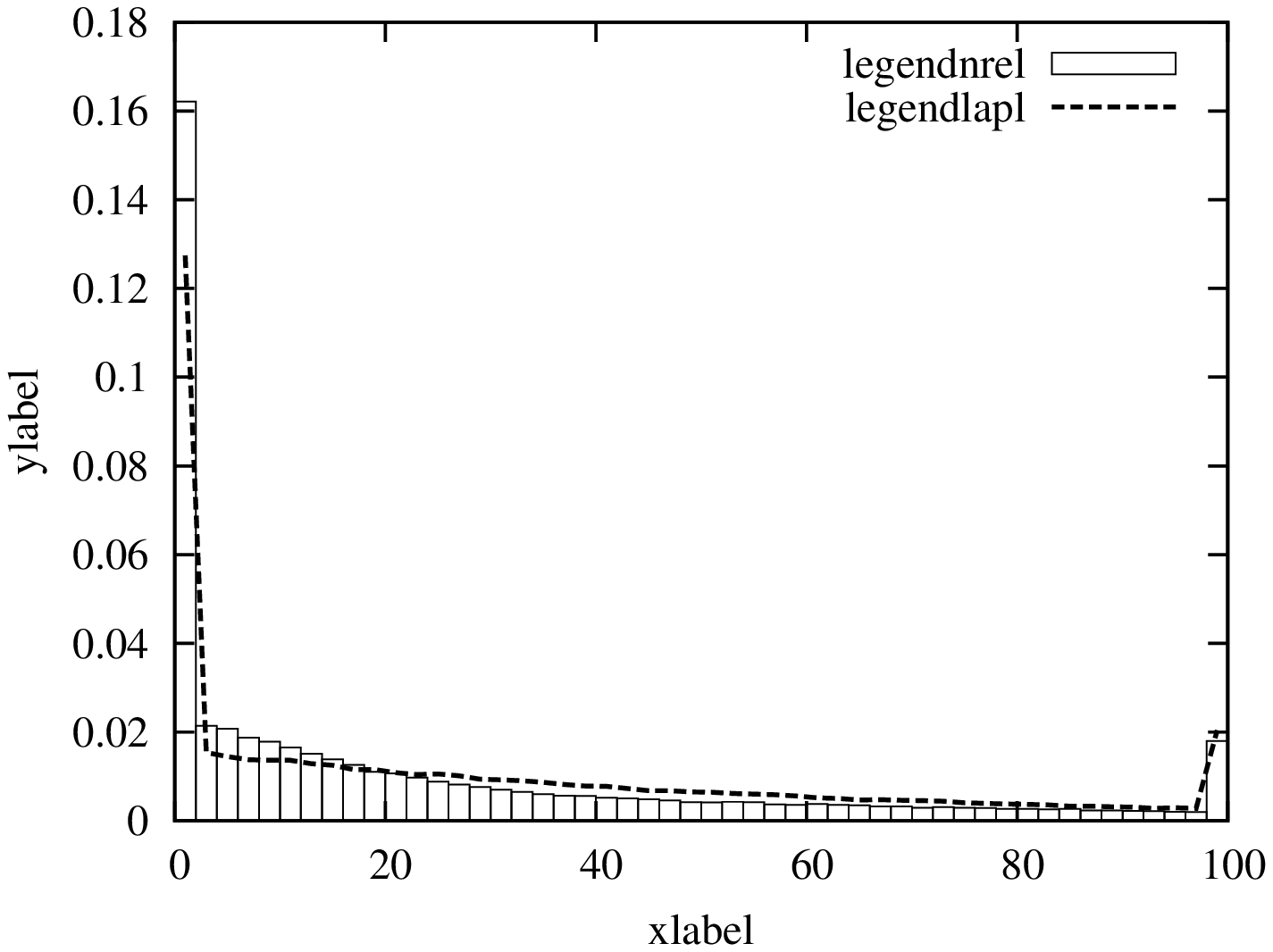}
\label{fig:fastramp-sdist}}
\subfloat[Fast-ramping generation]{
\footnotesize
\psfrag{xlabel}[c][c]{Fast-ramping generation (MW)}
\psfrag{ylabel}[c][c]{Probability density}
\psfrag{legendnrel}[r][r]{Wind dataset}
\psfrag{legendlapl}[r][r]{Laplace($\l$)}
\includegraphics[width=0.37\textwidth,angle=270]{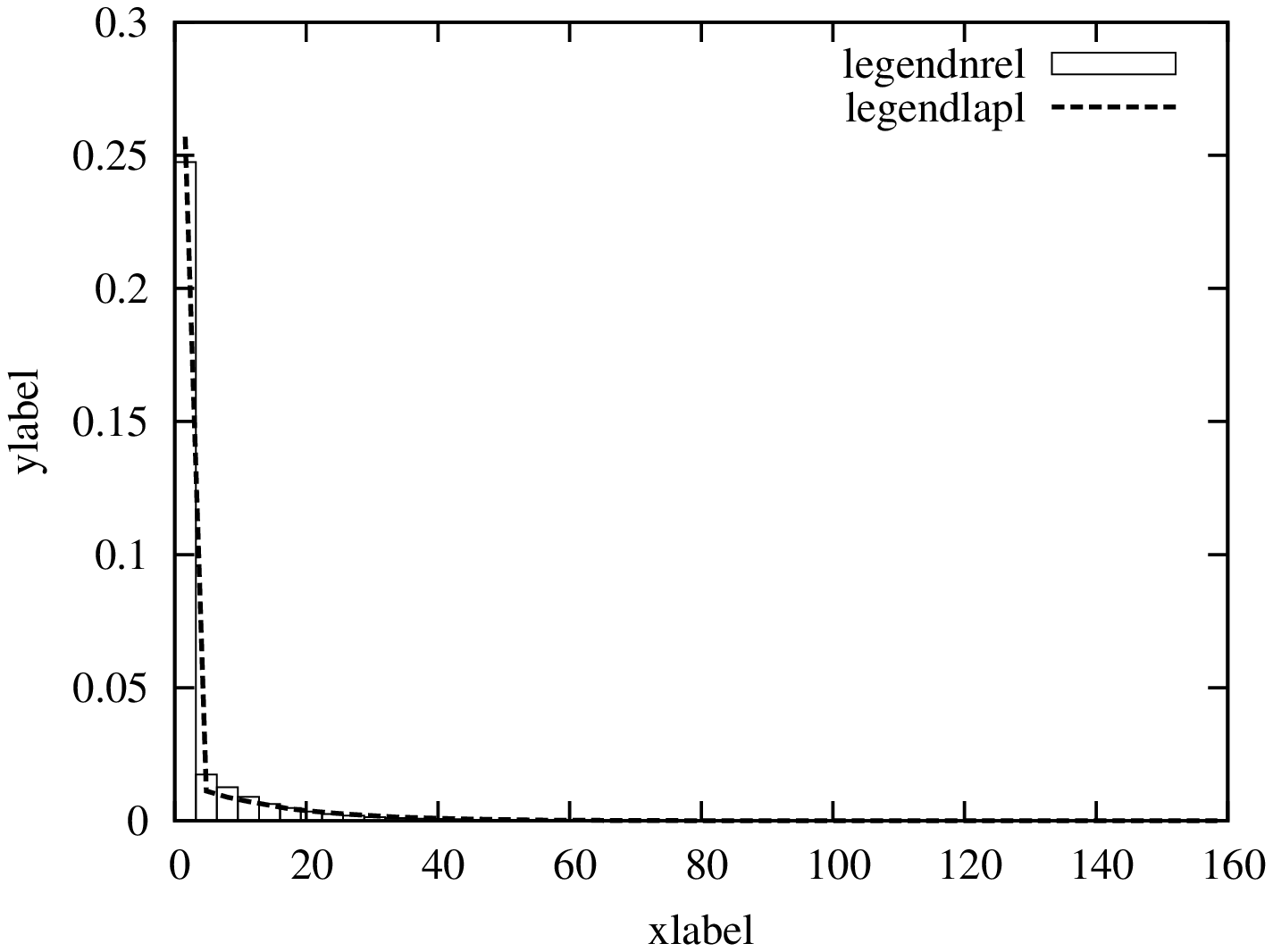}
\label{fig:fastramp-gdist}}
}
{Figures~\subref{fig:fastramp-sdist} and \subref{fig:fastramp-gdist} show the pdfs of the stored power and the fast-ramping generation, respectively, for $\a = 60\%$, $\Smax = 100$ MW, and $\Gmax = 160$ MW.
\label{fig:fastramp-dist}}
{}
\end{figure}

Figure~\ref{fig:fastramp-dist}\subref{fig:fastramp-sdist} compares the empirical pdf of the stored power of the simulated wind generation data to the stationary pdf under the Laplace distribution assumption in Proposition~\ref{prop:fastramp-dist}.
The corresponding empirical pdf of the fast-ramping generation and its stationary pdf are shown in Figure~\ref{fig:fastramp-dist}\subref{fig:fastramp-gdist}. Note again the simulation results corroborate well with the theory.

The loss of load probability depends on the tail of the cdf $F_\D$ of the net generation prediction error. However, the number of samples of the wind data is small, and thus it is difficult to compare the loss of load probabilities of the wind dataset and the Laplace-distributed prediction error sequence. To illustrate the loss of load probabilities for the residual power system, we will assume that the net generation prediction error sequence is distributed according to the corresponding best fit Laplace distribution. Then the expected loss of load probability can be expressed as
\begin{align*}
\Jcl(\pi,S_1) = \limsup_{n \to \infty} \frac{1}{n} \sum_{i=1}^n F_{\D}(- \Gmax - \ad S_i),
\end{align*}
where $F_{\D}$ is the cdf of the best fit Laplace distribution. 
Figure~\ref{fig:r}\subref{fig:lolp-r} compares the expected loss of load probabilities and using the three-year simulated wind data for various values of power storage capacities and round-trip storage efficiencies $\a=60\%$ and $80\%$.
For large $\Smax$, the average loss of load probability decreases exponentially in $\Smax$ as expected in Proposition~\ref{prop:lossofload-rate}.

% Two-threshold policy
%--------------------
\section{Two-threshold policy}
\label{sec:threshold}
In Section~\ref{sec:single-bus}, we showed that the policy $\pg$ minimizes the expected average conventional (or fast-ramping) generation---storage is discharged before conventional generation is used and conventional generation is never used to charge storage. We also showed that the policy $\pl$ minimizes the average loss of load probability---conventional generation is used before storage is discharged and conventional generation is used to keep storage as full as possible. In this section, we propose a general two-threshold policy, which includes $\pg$ and $\pl$ as extreme special cases. This policy tries to keep the stored power above a \emph{charging threshold} and tries not to discharge storage below a \emph{discharging threshold}.

Let $0 \leq \Sc \leq \Sd \leq \Smax$. 
The two-threshold policy is characterized by a charging threshold $\Sc$ and a discharging threshold $\Sd$ as given in Table~\ref{tab:threshold} and illustrated in Figure~\ref{fig:threshold} for $\Sc < \ac\Gmax$.
%The two-threshold policy is characterized by a charging threshold $\Sc$ and a discharging threshold $\Sd$ as given in Table~\ref{tab:threshold}.
When the prediction error $\D \geq 0$, the storage is charged as much as possible using excess net renewable generation. If the stored power after this charging is above $\Sc$, then conventional generation is not used. However, if it is below $\Sc$, then the storage is charged as close to $\Sc$ as possible using conventional generation. 
When the prediction error $\D < -\Gmax$, the storage must be discharged to balance the prediction error. If there is still unbalanced prediction error after the storage is discharged to $\Sd$, then fast-ramping generation is used such that the stored power is as close to $\Sd$ as possible. 
When the prediction error $-\Gmax \leq \D < 0$, the case where the stored power is either lower than $\Sc$ or higher than $\Sd$ is similar to the above cases. When the stored power is between $\Sc$ and $\Sd$, only the fast-ramping generation is used to balance the prediction error, and the stored power is unchanged.

Note that the optimal policies $\pg$ and $\pl$ for the two extreme cases in Theorem~\ref{thm:policy-generation} and~\ref{thm:policy-lossofload} are special cases of this policy with thresholds $(0,0)$ and $(\Smax,\Smax)$, respectively.

\begin{table*}[h]
\TABLE
{Two-threshold policy parameterized by $(\Sc,\Sd)$. \label{tab:threshold}}
{
\begin{minipage}{\textwidth}
\centering
\subfloat[$0 \leq S_i < \Sc$]{
\begin{tabular}{c|ccc|c}
$\pi_i$ & $G_i$ & $C_i$ & $D_i$ & $S_{i+1}$\\
\hline
$\displaystyle \frac{\Smax -S_i}{\ac} \leq \D_i$ & $0$ & $\displaystyle \frac{\Smax - S_i}{\ac}$ & $0$ & $\Smax$\\
$\displaystyle \frac{\Sc - S_i}{\ac} \leq \D_i < \frac{\Smax - S_i}{\ac}$ & $0$ & $\D_i$ & $0$ & $S_i + \ac \D_i$\\
$\displaystyle \frac{\Sc -S_i}{\ac} - \Gmax \leq \D_i < \frac{\Sc - S_i}{\ac}$ & $\displaystyle \frac{\Sc -S_i}{\ac} - \D_i$ & $\displaystyle \frac{\Sc - S_i}{\ac}$ & $0$ & $\Sc$\\
$\displaystyle -\Gmax \leq \D_i < \frac{\Sc - S_i}{\ac} - \Gmax$ & $\Gmax$ & $\Gmax + \D_i$ & $0$ & $S_i + \ac(\Gmax+\D_i)$\\
$-\Gmax - \ad S_i \leq \D_i < -\Gmax$ & $\Gmax$ & $0$ & $- \Gmax - \D_i$ & $\displaystyle S_i + \frac{\Gmax+\D_i}{\ad}$\\
$\D_i < -\Gmax - \ad S_i$ & $\Gmax$ & $0$ & $\ad S_i$ & $0$\\
\hline
\end{tabular}
\label{tab:threshold-low}} \\
\subfloat[$\Sc \leq S_i \leq \Sd$]{
\begin{tabular}{c|ccc|c}
$\pi_i$ & $G_i$ & $C_i$ & $D_i$ & $S_{i+1}$\\
\hline
$\displaystyle \frac{\Smax -S_i}{\ac} \leq \D_i$ & $0$ & $\displaystyle \frac{\Smax - S_i}{\ac}$ & $0$ & $\Smax$ \\
$\displaystyle 0 \leq \D_i < \frac{\Smax - S_i}{\ac}$ & $0$ & $\D_i$ & $0$ & $S_i + \ac \D_i$ \\
$-\Gmax \leq \D_i < 0$ & $-\D_i$ & $0$ & $0$ & $S_i$ \\
$-\Gmax - \ad S_i \leq \D_i < -\Gmax$ & $\Gmax$ & $0$ & $- \Gmax - \D_i$ & $\displaystyle S_i + \frac{\Gmax+\D_i}{\ad}$ \\
$\d < -\Gmax - \ad S_i$ & $\Gmax$ & $0$ & $\ad S_i$ & $0$ \\
\hline
\end{tabular}
\label{tab:threshold-med}} \\
\subfloat[$\Sd < S_i \leq \Smax$]{
\begin{tabular}{c|ccc|c}
$\pi_i$ & $G_i$ & $C_i$ & $D_i$ & $S_{i+1}$\\
\hline
$\displaystyle \frac{\Smax -S_i}{\ac} \leq \D_i$ & $0$ & $\displaystyle \frac{\Smax - S_i}{\ac}$ & $0$ & $\Smax$ \\
$\displaystyle 0 \leq \D_i < \frac{\Smax - S_i}{\ac}$ & $0$ & $\D_i$ & $0$ & $S_i + \ac \D_i$ \\
$-\ad(S_i - \Sd) \leq \D_i < 0$ & $0$ & $0$ & $-\D_i$ & $S_i + \D_i/\ad$ \\
$-\ad(S_i-\Sd) - \Gmax \leq \D_i < -\ad(S_i - \Sd)$ & $-\D_i - \ad(S_i - \Sd)$ & $0$ & $\ad(S_i - \Sd)$ & $\Sd$ \\
$-\Gmax - \ad S_i \leq \D_i < -\ad(S_i - \Sd) -\Gmax$ & $\Gmax$ & $0$ & $- \Gmax - \D_i$ & $\displaystyle S_i + \frac{\Gmax+\D_i}{\ad}$ \\
$\D_i < -\Gmax - \ad S_i$ & $\Gmax$ & $0$ & $\ad S_i$ & $0$ \\
\hline
\end{tabular}
\label{tab:threshold-high}}
\end{minipage}
}
{}
\end{table*}

\begin{figure}[h]
\FIGURE
{\begin{tikzpicture}[scale=0.6,font=\footnotesize]
\node (s) at (8.8,0) {$S_i$};
\node (d) at (0,5.8) {$\D_i$};
\draw[-latex] (0,-5.5) -- (0,5.5);
\draw[-latex] (-0.5,0) -- (8.5,0);
\draw (8,5.5) -- (8,-5.5);
\draw (0,5) -- (8,0);
\node[anchor=east] at (0,5) {$\Smax/\ac$};
\node at (6,3) {$\Smax$};
\node at (3,1.5) {$S_i + \ac \D_i$};
\draw (0,1.25) -- (2,0);
\node[anchor=east] at (0,1.25) {$\Sc/\ac$};
\node at (0.6,0.3) {$\Sc$};
\node at (1,-0.5) {$\Sc$};
\node at (3,-1) {$S_i$};
\draw (0,-0.75) -- (2,-2);
\node[anchor=east] at (0,-0.65) {$\Sc/\ac-\Gmax$};
\node[anchor=east] at (-0.4,-1.5) {$S_i +\ac(\Gmax+\D_i)$}; 
\draw[-latex] (-0.5,-1.5) -- (0.5,-1.5);
\node[anchor=west] at (8.4,-0.7) {$S_i +\D_i/\ad$};
\draw[-latex] (8.5,-0.7) -- (7,-0.7);
\node[anchor=west] at (8.4,-4.4) {$S_i + (\Gmax+\D_i)/\ad$};
\draw[-latex] (8.5,-4.4) -- (7,-4.4);
\draw (0,-2) -- (4,-2);
\node[anchor=east] at (0,-2.1) {$-\Gmax$};
\draw (0,-2) -- (8,-5.2);
\node[anchor=west] at (8,-5.3) {$- \Gmax - \ad \Smax$};
\node at (3,-4.5) {$0$};
\draw (2,-2) -- (2,0);
\draw (4,0) -- (4,-2);
\draw (4,0) -- (8,-1.6);
\node[anchor=west] at (8,-1.7) {$-\ad(\Smax-\Sd)$};
\node at (6,-1.7) {$\Sd$};
\draw (4,-2) -- (8,-3.6);
\node[anchor=west] at (8,-3.5) {$-\ad(\Smax-\Sd) - \Gmax$};
\end{tikzpicture}}
{Illustration of the two-threshold policy.
\label{fig:threshold}}
{}
\end{figure}

Next we consider the general weighted-sum average cost:
\begin{align*}
\Jc(\pi,S_1) = \rho_1 \Jcg(\pi,S_1) + \rho_2 \Jcl(\pi,S_1),
\end{align*}
where $\rho_1 \geq 0$ and $\rho_2 \geq 0$. 
In the following proposition, we show that the two-threshold policy is optimal for the residual power system dynamic program with the above  general weighted-sum cost and two slots. 

\begin{proposition}
\label{prop:threshold}
If the pdf of the prediction error $f_\D$ increases on $(-\infty,0]$, then there exist ${\Sc}_i$ and ${\Sd}_i$ such that $0 \leq {\Sc}_i \leq {\Sd}_i \leq \Smax$, and the optimal policy for the general residual power system dynamic program with two slots is the two-threshold policy with parameters $({\Sc}_i,{\Sd}_i)$ for time $i = 1, 2$.
\end{proposition}
\medskip
The proof of this proposition is given in Appendix~\ref{app:threshold}.

% numerical result
%--------------------
To demonstrate the two-threshold policy, we implemented a dynamic programming method by discretizing the state space and then running the value iteration
%%%(Bertsekas 2007)
\citep{Bertsekas2007}. 
For the values of $\Gmax$, $\Smax$, $\rho_0$, and $\rho_1$ used in the following numerical examples, we find that the policy obtained from the value iteration is a discretized two-threshold policy. 

\begin{figure}[h]
\FIGURE
{\footnotesize
\psfrag{xlabel}[c][c]{Loss of load probability $\Jcl$}
\psfrag{ylabel}[c][c]{Fast-ramping generation $\Jcg$ (MW)}
\psfrag{legendnrel1}[r][r]{Wind dataset $\Smax = 0$}
\psfrag{legendnrel2}[r][r]{Wind dataset $\Smax = 50$}
\psfrag{legendnrel3}[r][r]{Wind dataset $\Smax = 100$}
\includegraphics[width=0.45\textwidth,angle=270]{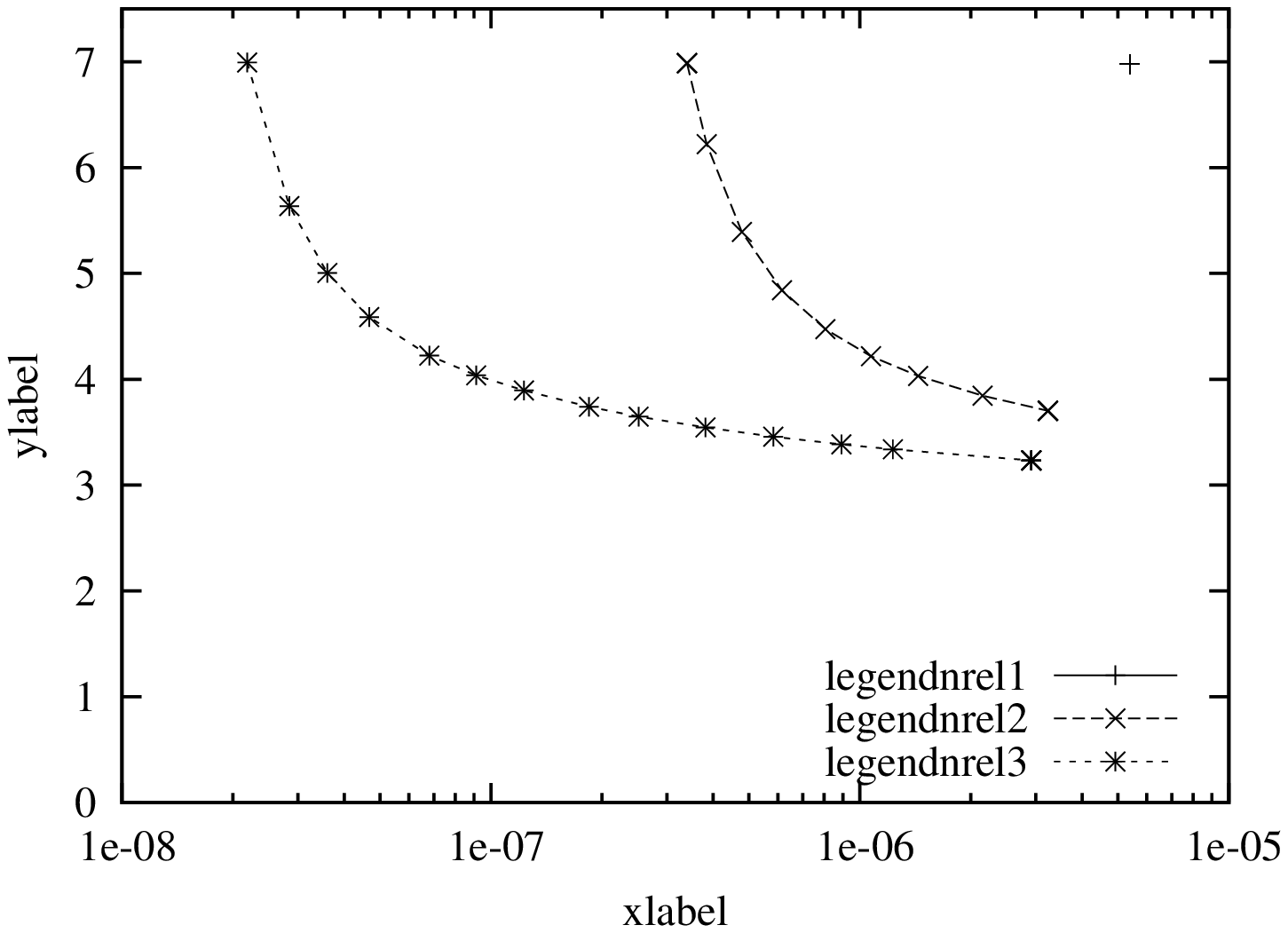}}
{The tradeoff between the fast-ramping generation and the loss of load probability for $\Gmax = 160$ MW under two-threshold policies.
\label{fig:pareto}}
{}
\end{figure}

\noindent{\em Tradeoff between $\Jcg$ and $\Jcl$}: Figure~\ref{fig:pareto} shows the tradeoff between the fast-ramping generation and the loss of load probability for no storage and for storage capacities $\Smax = 50$ MW and $\Smax =100$ MW with fast-ramping generation capacity $\Gmax = 160$ MW. Note that the results for the Laplace pdf corroborate very well with the simulated wind generation data. As shown in the figure the loss of load probability is improved by more than two orders of magnitude by using power storage capacity less than 5 standard deviations of the prediction error.

\begin{comment}
\noindent{\em Empirical distribution of stored power}: The empirical pdf of the stored power for $(\rho_0,\rho_1) = (10^{-6}, 1)$ is shown in Figure~\ref{fig:distri-storage}. 
From Figure~\ref{fig:threshold}, we expect the probability concentrates at $0$, $\Sc$, $\Sd$, and $\Smax$. The empirical pdf in Figure~\ref{fig:distri-storage} indeed has two peaks at $\Sd$ and $\Smax$. The probability that the stored power is below $\Sd$ is very small since the transition probability of the stored power from a value above $\Sd$ to below $\Sd$ is less than $F_\D(-\Gmax) = 5.4\cdot 10^{-6}$.
\end{comment}

\begin{figure}[h]
\FIGURE
{\footnotesize
\psfrag{xlabel}[c][c]{Power storage capacity $\Smax$ (MW)}
\psfrag{ylabel}[c][c]{Fast-ramping power capacity $\Gmax$ (MW)}
\psfrag{legendnrel}[r][r]{Wind dataset}
\psfrag{legendlapl}[r][r]{Laplace($\l$)}
\includegraphics[width=0.45\textwidth,angle=270]{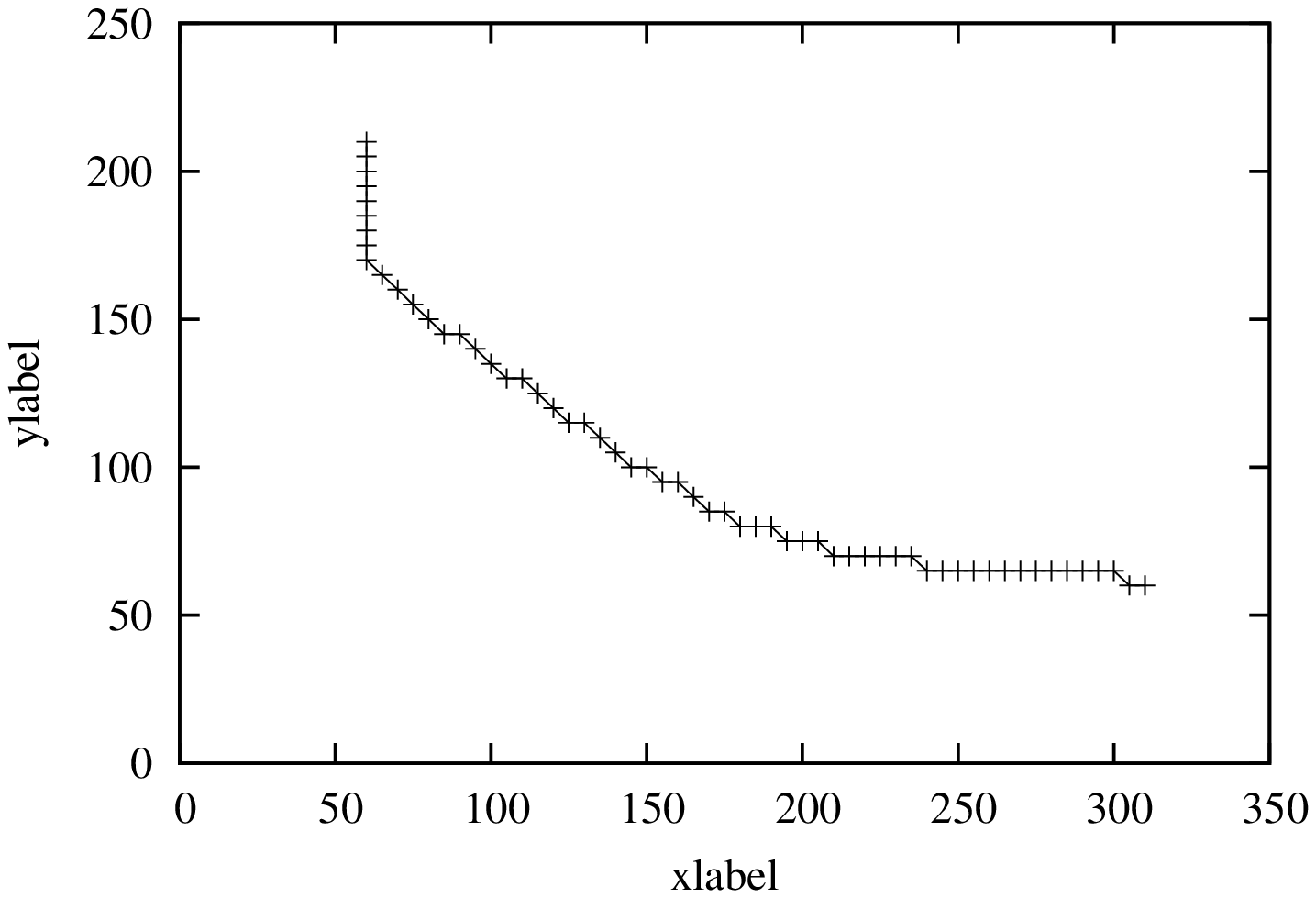}}
{The tradeoff between the fast-ramping generation power capacity and the power storage capacity for $\Jcg \leq 3.6$ MW and $\Jcl \leq 2\cdot 10^{-6}$.
\label{fig:plan}}
{}
\end{figure}

\noindent{\em Tradeoff between $\Smax$ and $\Gmax$}: In Figure~\ref{fig:plan}, we compare the two ways of mitigating renewable energy variability; using fast-ramping generation and using fast-response storage. We fix the expected average fast-ramping generation at $3.6$ MW (corresponding to 80\% maximum reduction in the fast-ramping generation) and the loss of load probability at $2\cdot 10^{-6}$ (corresponds to one loss of load event every 10 years). To achieve these goals with minimum power storage capacity, we need $\Gmax = 170$ MW and $\Smax = 60$ MW. To reduce the fossil fuel generation and to achieve the same goals, we can replace 1 MW of $\Gmax$ with 1.3 MW of $\Smax$.

\begin{comment}
\begin{figure}[h]
\FIGURE
{\footnotesize
\psfrag{xlabel}[c][c]{Stored power (MW)}
\psfrag{ylabel}[c][c]{Probability density}
\psfrag{legendnrel}[r][r]{Wind dataset}
\psfrag{legendlapl}[r][r]{Laplace($\l$)}
\includegraphics[width=0.3\textwidth,angle=270]{distri-storage}}
{The pdfs of the stored power under the two-threshold policy with $\a = 60\%$, $\Gmax = 160$ MW, $\Smax = 100$ MW, $\Sc = 12$ MW, and $\Sd = 38$ MW.
\label{fig:distri-storage}}
{}
\end{figure}
\end{comment}

\section{Constrained $\Cmax$ and $\Dmax$}
\label{sec:constrained}
% constrained \Cmax and \Dmax
%--------------------
In previous sections, we assumed the  rated storage power conversion $\Cmax$ and the rated storage output power $\Dmax$ to be unconstrained. In this section, we relax this condition and show that Theorems~\ref{thm:policy-generation} and~\ref{thm:policy-lossofload} can be extended to find the optimal policies with constraints on $\Cmax$ and  $\Dmax$.  

Consider stochastic program I with the additional constraints
\begin{align}
\Cmax & \leq \frac{1}{\ac} \Smax, \qquad \Dmax \leq \ad \Smax. \label{equ:constr-cd}
\end{align}

In this following, we show that the optimal policy is a slight generalization of the policy for the unconstrained case in Theorem~\ref{thm:policy-generation}. %

\begin{comment}
\begin{theorem}
\label{thm:policy-generation-constr}
The optimal policy $\pg$ for constrained $\Cmax$ and $\Dmax$ is given in Table~\ref{tab:policy-generation-constr}.
\end{theorem}
\medskip
\end{comment}

\begin{theorem}
\label{thm:policy-generation-constr}
The optimal policy $\pi^{(3)}$ for stochastic program I with constrained $\Cmax$ and $\Dmax$ is equal to
\begin{align*}
C_i^{(3)} & = \min\{ \Cg_i, \Cmax \}, \qquad
D_i^{(3)} = \min\{ \Dg_i, \Dmax \}, \qquad
G_i^{(3)} = \min\{ [C_i^{(3)} - D_i^{(3)} - \D_i]^+, \Gmax \}.
\end{align*}
\end{theorem}
\medskip
The proof of this theorem is given in Appendix~\ref{app:single-bus}. 
Note that the above optimal policy reduces to the policy in Theorem~\ref{thm:policy-generation} when $\ac\Cmax = \Dmax/\ad = \Smax$. 

\begin{comment}
\begin{table}[h]
\TABLE
{Optimal policy in Theorem~\ref{thm:policy-generation-constr} and corresponding stored power for constrained $\Cmax$ and $\Dmax$.
\label{tab:policy-generation-constr}}
{\begin{tabular*}{\textwidth}{@{\extracolsep{\fill}}c|ccc|c}
$\pg$ & $\Gg_i$ & $\Cg_i$ & $\Dg_i$ & $S_{i+1}$ \\
\hline
$\displaystyle \frac{\Smax-S_i}{\ac} \leq \d_i$, $\Smax-\ac\Cmax \leq S_i$ & $0$ & $\displaystyle \frac{\Smax-S_i}{\ac}$ & $0$ & $\Smax$ \\
$\Cmax \leq \D_i, S_i < \Smax - \ac\Cmax$ & $0$ & $\Cmax$ & $0$ & $S_i + \ac \Cmax$ \\
$\displaystyle 0 \leq \D_i < \min\left\{\frac{\Smax-S_i}{\ac},\Cmax\right\}$ & $0$ & $\D_i$ & $0$ & $S_i + \ac \D_i$\\
$\max\{-\ad S_i,-\Dmax\} \leq \D_i < 0$ & $0$ & $0$ & $-\D_i$ & $S_i + \D_i/\ad$ \\
$- \Gmax - \ad S_i \leq \D_i < -\ad S_i$, $S_i \leq \Dmax/\ad$ & $-\D_i - \ad S_i$ & $0$ & $\ad S_i$ & $0$ \\
$- \Gmax - \Dmax \leq \D_i < -\Dmax$, $\Dmax/\ad < S_i$ & $-\D_i - \Dmax$ & $0$ & $\Dmax$ & $S_i - \Dmax/\ad$ \\
$\D_i < - \Gmax - \ad S_i$, $S_i \leq \Dmax/\ad$ & $\Gmax$ & $0$ & $\ad S_i$ & $0$\\
$\D_i < - \Gmax - \Dmax$, $\Dmax/\ad < S_i$ & $\Gmax$ & $0$ & $\Dmax$ & $S_i - \Dmax/\ad$\\
\hline
\end{tabular*}}
{}
\end{table}
\end{comment}

Similarly, we can consider stochastic program II with constraints in~\eqref{equ:constr-cd} and obtain a generalization of the policy for the unconstrained case in 
Theorem~\ref{thm:policy-lossofload}.

\begin{comment}
\begin{theorem}
\label{thm:policy-lossofload-constr}
For constrained $\Cmax$ and $\Dmax$, the optimal policy $\pl$ is given in Table~\ref{tab:policy-lossofload-constr}, and it is stationary.
\end{theorem}
\medskip
\end{comment}

\begin{theorem}
\label{thm:policy-lossofload-constr}
The optimal policy $\pi^{(4)}$ for stochastic program II with constrained $\Cmax$ and $\Dmax$ is equal to
\begin{align*}
C_i^{(4)} & = \min\{ \Cl_i, \Cmax \}, \qquad
D_i^{(4)} = \min\{ \Dl_i, \Dmax \}, \qquad
G_i^{(4)} = \min\{ [C_i^{(4)} - D_i^{(4)} - \D_i]^+, \Gmax \}.
\end{align*}
\end{theorem}
\medskip
The proof of this theorem is given in Appendix~\ref{app:single-bus}.

\begin{comment}
\begin{table}[t]
\TABLE
{Optimal policy in Theorem~\ref{thm:policy-lossofload-constr} and corresponding stored power for constrained $\Cmax$ and $\Dmax$.
\label{tab:policy-lossofload-constr}}
{\begin{tabular*}{\textwidth}{@{\extracolsep{\fill}}c|ccc|c}
$\pl$ & $\Gl_i$ & $\Cl_i$ & $\Dl_i$ & $S_{i+1}$ \\
\hline
$\displaystyle \frac{\Smax-S_i}{\ac} - \Gmax \leq \D_i$, $\Smax-\ac\Cmax \leq S_i$ & $\displaystyle \left(\frac{\Smax-S_i}{\ac}-\D_i \right)^+$ & $\displaystyle \frac{\Smax-S_i}{\ac}$ & $0$ & $\Smax$ \\
$\Cmax - \Gmax \leq \D_i$, $s < \Smax - \ac\Cmax$ & $(\Cmax-\D_i)^+$ & $\Cmax$ & $0$ & $S_i + \ac \Cmax$ \\
$\displaystyle -\Gmax \leq \D_i < \min\left\{\frac{\Smax-S_i}{\ac},\Cmax\right\} - \Gmax$ & $\Gmax$ & $\D_i + \Gmax$ & $0$ & $S_i + \ac(\D_i+\Gmax)$ \\
$\max\{-\ad S_i,-\Dmax\} - \Gmax \leq \D_i < -\Gmax$ & $\Gmax$ & $0$ & $-\D_i - \Gmax$ & $\displaystyle S_i + \frac{\D_i+\Gmax}{\ad}$ \\
$\D_i < - \Gmax - \ad S_i$, $S_i \leq \Dmax/\ad$ & $\Gmax$ & $0$ & $\ad S_i$ & $0$ \\
$\D_i < - \Gmax - \Dmax$, $\Dmax/\ad < S_i$ & $\Gmax$ & $0$ & $\Dmax$ & $S_i - \Dmax/\ad$ \\
\hline
\end{tabular*}}
{}
\end{table}
\end{comment}

\begin{figure}[h]
\FIGURE
{
\subfloat[]{
\footnotesize
\psfrag{xlabel}[c][c]{$\ac\Cmax/\Smax = \Dmax/\ad\Smax$}
\psfrag{ylabel}[c][c]{Conventional generation $\Jcg$ (MW)}
\psfrag{legend-nrel6}[r][r]{$\a = 60\%$}
\psfrag{legend-nrel8}[r][r]{$\a = 80\%$}
\includegraphics[width=0.37\textwidth,angle=270]{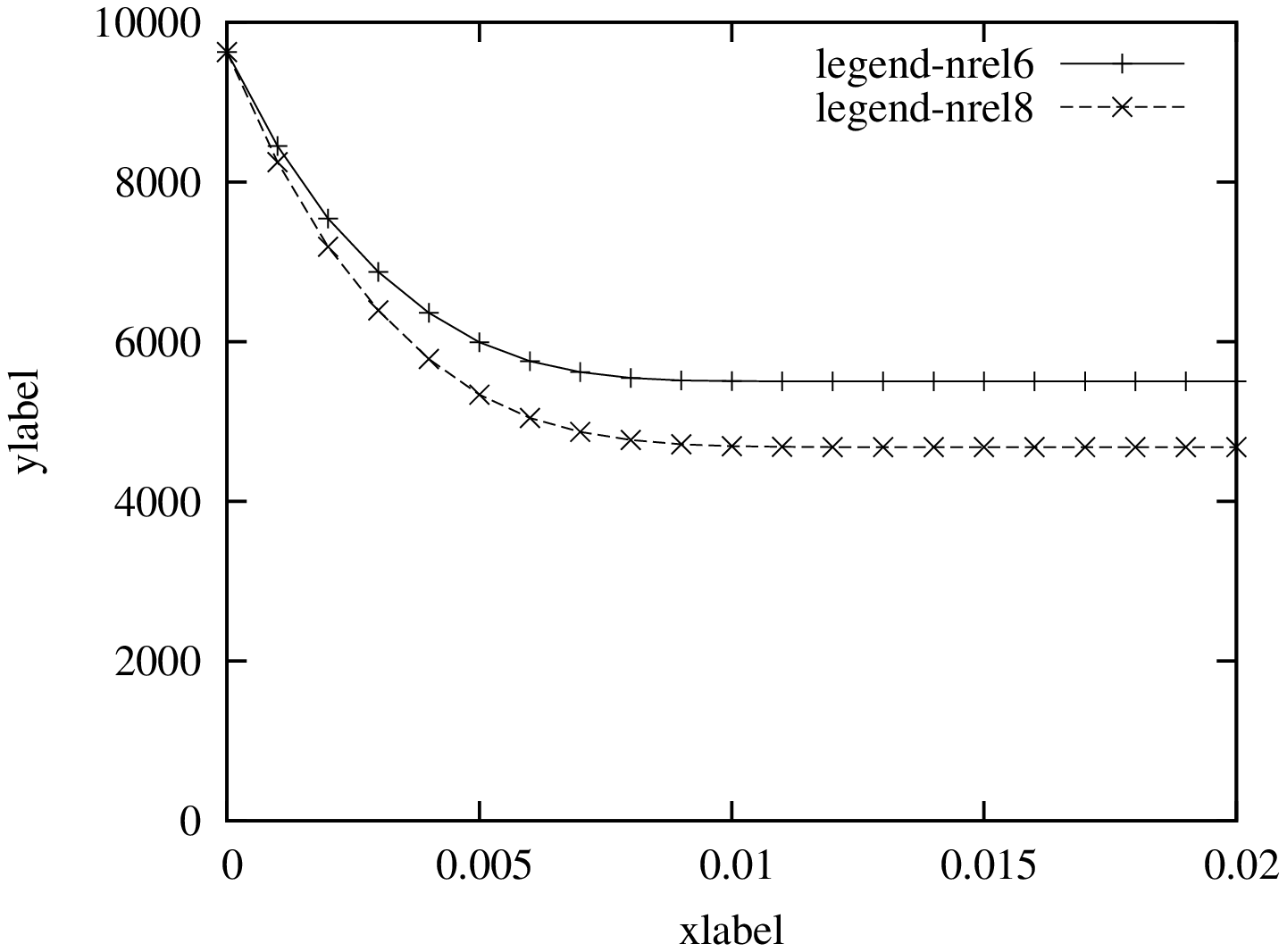}
\label{fig:generation-constr-s}
}
\subfloat[]{
\footnotesize
\psfrag{xlabel}[c][c]{$\ac\Cmax/\Smax = \Dmax/\ad\Smax$}
\psfrag{ylabel}[c][c]{Loss of load probability $\Jcl$}
\psfrag{legend-nrel6}[r][r]{$\a = 60\%$}
\psfrag{legend-nrel8}[r][r]{$\a = 80\%$}
\includegraphics[width=0.37\textwidth,angle=270]{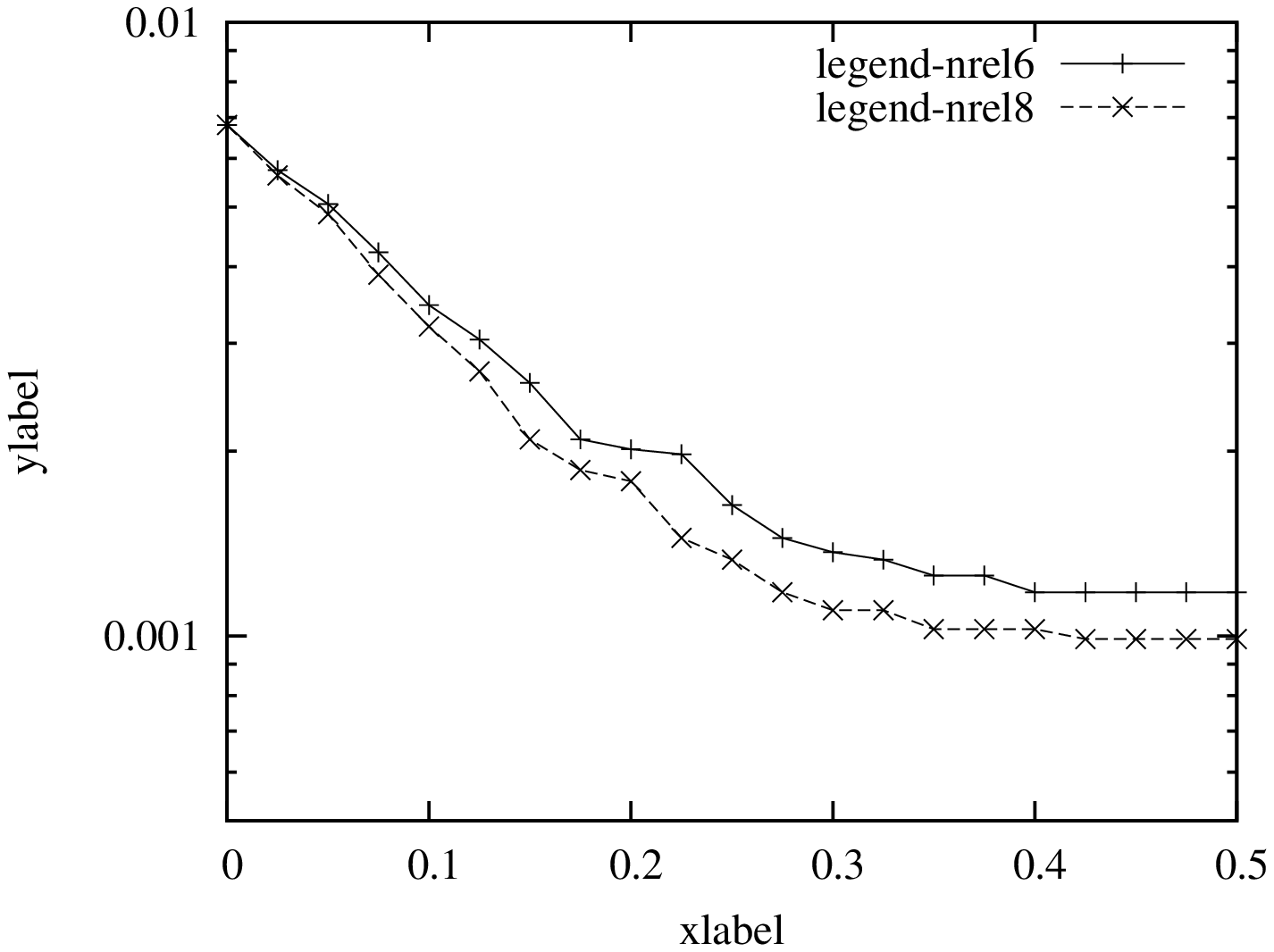}
\label{fig:lolp-constr-s}
}
}
{Figures~\subref{fig:generation-constr-s} and \subref{fig:lolp-constr-s} show the minimum expected average conventional generation and the minimum expected loss of load probability, respectively, versus power storage capacity for round-trip efficiencies $\a=60\%$ and $80\%$ and $\Cmax=\Gmax =160$ MW.
\label{fig:constr-s}}
{}
\end{figure}

%\subsection{Numerical results} 
%--------------------
Figure~\ref{fig:constr-s}\subref{fig:generation-constr-s} plots the minimum average conventional generation versus $\ac \Cmax/\Smax = \Dmax/ \ad \Smax$ for the NREL wind dataset and the CAISO load dataset. Note that most of the reduction in conventional generation in the unconstrained case can be achieved by $\Cmax$ and $\Dmax$ less than $5\%$ of $\Smax$.

Figure~\ref{fig:constr-s}\subref{fig:lolp-constr-s} plots the minimum expected average loss of load probability versus $\ac \Cmax/\Smax = \Dmax/ \ad \Smax$. Note that loss of load may occur even when the stored power is high due to constrained $\Dmax$. Thus, the loss of load probability can be much higher for small $\Cmax$ and $\Dmax$ than for the unconstrained case.

% average convetional generation v.s. Cmax

\section{Conclusion}

This paper aimed to answer questions concerning the limits on the benefit of energy storage to renewable energy: How much can storage help? how much storage is needed? and what are the optimal control strategies for storage? To answer these questions, we formulated a single-bus power system with storage as stochastic program and established the optimal policies when the cost function is the expected average conventional generation and when it is the expected average loss of load probability. We proposed a general two-threshold policy for which these optimal policies are extreme special cases. We obtained refined analytical results by considering the multi-timescale operation of the grid. The results lead to the following potentially useful conclusions: 
\begin{itemize}
\item Using energy storage reduces fast-ramping generation by a factor up to the reciprocal of the round-trip inefficiency.
\item 80\% of the reduction in fast-ramping generation can be achieved using storage capacity equal to only four standard deviations of the net renewable generation prediction error.
\item The loss of load probability can be reduced by an unbounded factor as the energy storage capacity increases.
\item The loss of load probability can be reduced by an order of magnitude using storage capacity equal to only two standard deviations of the prediction error. 
\item The Laplace distribution, which makes the analysis far more tractable, appears to be a reasonable approximation of the short time-scale prediction error of wind energy generation.
\item Most reduction in conventional generation for unconstrained charging and discharging rates can be achieved by relatively small such rates.
\end{itemize}
We corroborated our assumptions and analytical results with simulated wind data.

There are many open questions suggested by this work:
What is the minimum expected average fast-ramping generation for constrained $\Cmax$ and $\Dmax$ under the Laplace assumption?
What is the minimum expected average loss of load probability for limited $\Gmax$ and $\Smax$?
What is the natural cost function for which the two-threshold policy is optimal?
What is the optimal policy when the ramping constraints of the fast-ramping generation are considered?

% Appendix here
% Options are (1) APPENDIX (with or without general title) or
%             (2) APPENDICES (if it has more than one unrelated sections)
% Outcomment the appropriate case if necessary
%
% \begin{APPENDIX}{<Title of the Appendix>}
% \end{APPENDIX}
%
%   or
%
% \begin{APPENDICES}
% \section{<Title of Section A>}
% \section{<Title of Section B>}
% etc
% \end{APPENDICES}

\begin{APPENDICES}
\section{Proofs for the single-bus power system results}
\label{app:single-bus}

% proof of optimal policy i
%--------------------
\proof{Proof of Theorems~\ref{thm:policy-generation} and~\ref{thm:policy-generation-constr}.}
Since Theorem~\ref{thm:policy-generation} is a special case of Theorem~\ref{thm:policy-generation-constr}, we only need to prove the general case. Note that the proof of Lemma~\ref{lm:jg-upper} holds for constrained $\Cmax$ and $\Dmax$ since the policies in~\eqref{equ:jg-upper-policy} and~\eqref{equ:jg-lower-policy} have smaller charging and discharging power than the optimal policies.
Thus, we only need to show that for any policy $\pi_{k-1}$,
\begin{align*}
J_{k-1} & = G_{k-1} + \Jg_k\left(\pi^*,S_{k-1} + \ac C_{k-1} - \frac{1}{\ad} D_{k-1}  \right) \\
& \geq  \Gg_{k-1} + \Jg_k\left(\pi^*,S_{k-1} + \ac \Cg_{k-1} - \frac{1}{\ad} \Dg_{k-1} \right) = \Jg_{k-1}. 
\end{align*}
Now we consider the following seven cases.
\begin{enumerate}
\item $(\Smax-S_{k-1})/\ac \leq \D_{k-1}$ and $\Smax - \ac\Cmax \leq S_{k-1}$. By Lemma~\ref{lm:jg-upper}, 
\begin{align*}
J_{k-1} \geq \Jg_k(\pi^*,\Smax) = \Jg_{k-1}.
\end{align*}

\item $\Cmax \leq \D_{k-1}$ and $S_{k-1} < \Smax - \ac\Cmax$. By Lemma~\ref{lm:jg-upper},
\begin{align*}
J_{k-1} \geq \Jg_k(\pi^*, S_{k-1}+\ac\Cmax) = \Jg_{k-1}.
\end{align*}

\item $0 \leq \D_{k-1} \leq \min\{(\Smax-S_{k-1})/\ac, \Cmax\}$.
\begin{itemize}
\item If $\a C_{k-1} - D_{k-1} \leq \a \Cg_{k-1} - \Dg_{k-1}$, then by Lemma~\ref{lm:jg-upper},
\begin{align*}
J_{k-1} \geq \Jg_k\left(\pi^*,S_{k-1} + \ac \Cg_{k-1} - \frac{1}{\ad} \Dg_{k-1} \right) = \Jg_{k-1}.
\end{align*}

\item If $\a C_{k-1} - D_{k-1} > \a \Cg_{k-1} - \Dg_{k-1}$, then by Lemma~\ref{lm:jg-upper},
\begin{align*}
J_{k-1} - \Jg_{k-1} & \geq C_{k-1} - \D_{k-1} - D_{k-1} \\
& \qquad - \ad \left( \ac (C_{k-1} - \D_{k-1}) - \frac{1}{\ad} D_{k-1} \right)  \geq 0.
\end{align*}
\end{itemize}

\item $- \min\{\ad S_{k-1},\Dmax\} \leq \D_{k-1} < 0$. The proof is the same as the previous case.

\item $-\Gmax - \ad S_{k-1} \leq \D_{k-1} < -\ad S_{k-1}$ and $S_{k-1} \leq \Dmax/\ad$. By Lemma~\ref{lm:jg-upper},
\begin{align*}
J_{k-1} - \Jg_{k-1} & \geq C_{k-1} - D_{k-1} + \ad S_{k-1} \\
& \qquad - \ad \left( S_{k-1} + \ac C_{k-1} - \frac{1}{\ad} D_{k-1} \right) \geq 0.
\end{align*}

\item $-\Gmax - \Dmax \leq \D_{k-1} < -\Dmax$ and $\Dmax/\ad < S_{k-1}$.
\begin{align*}
J_{k-1} - \Jg_{k-1} & \geq C_{k-1} - D_{k-1} + \Dmax \\
& \qquad - \ad \left( \ac C_{k-1} - \frac{1}{\ad}(D_{k-1}-\Dmax) \right) \geq 0.
\end{align*}

\item $\D_{k-1} < -\Gmax - \min\{\ad S_{k-1}, \Dmax\}$. By the definition of loss of load, $J_{k-1} = \Jg_{k-1}$.
\end{enumerate}

By induction we show that the stationary policy $\pg$ is optimal for the finite horizon stochastic program. Then for any $n \in \N$ and policy $\pi$, we have
\begin{align*}
\frac{1}{n} \Jg_1(\pg, S_1) \leq \frac{1}{n} \Jg_1(\pi, S_1),
\end{align*} 
and thus
\begin{align*}
\Jg(\pg, S_1) = \limsup_{n\to\infty} \frac{1}{n} \Jg_1(\pg,S_1) 
\leq \limsup_{n\to\infty} \frac{1}{n} \Jg_1(\pi,S_1) = \Jg(\pi,S_1).
\end{align*}

\endproof

% proof of optimal policy ii
%--------------------
\proof{Proof of Theorems~\ref{thm:policy-lossofload} and~\ref{thm:policy-lossofload-constr}.}
Since Theorem~\ref{thm:policy-lossofload} is a special case of Theorem~\ref{thm:policy-lossofload-constr}, we only need to prove the general case. Note that the property of the cost-to-go function in~\eqref{equ:jl-upper} holds for constrained $\Cmax$ and $\Dmax$ since the policy in~\eqref{equ:jg-upper-policy} have smaller charging and discharging power than the optimal policy.
Thus, we need to show that for any policy $\pi_{k-1}$,
\begin{align*}
J_{k-1} & = 1_{\{\D_{k-1} < -\Gmax - \ad S_{k-1}\}} + \Jl_k\left(\pi^*, S_{k-1} + \ac C_{k-1} - \frac{1}{\ad}D_{k-1} \right) \\
& \geq 1_{\{\D_{k-1} < -\Gmax - \ad S_{k-1}\}} + \Jl_k\left( \pi^*, S_{k-1} + \ac \Cl_{k-1} - \frac{1}{\ad}\Dl_{k-1} \right) 
= \Jl_{k-1}.
\end{align*}
Since the first terms in $J_{k-1}$ and $\Jl_{k-1}$ are equal, by~\eqref{equ:jl-upper}, we only need to show that
\begin{align*}
S_k & = S_{k-1} + \ac C_{k-1} - \frac{1}{\ad}D_{k-1} \\
& \geq S_{k-1} + \ac \Cl_{k-1} - \frac{1}{\ad}\Dl_{k-1} = \Sl_k.
\end{align*}

Now we consider the following seven cases.
\begin{enumerate}
\item If $(\Smax - S_{k-1})/\ac - \Gmax \leq \D_{k-1}$ and $\Smax - \ac\Cmax \leq S_{k-1}$, then $S_k \leq \Smax = \Sl_k$.

\item If $\Cmax - \Gmax \leq \D_{k-1}$ and $S_{k-1} < \Smax - \ac\Cmax$, then $S_k \leq S_{k-1} + \ac\Cmax = \Sl_k$.

\item If $-\Gmax \leq \D_{k-1} < \min\{ (\Smax - S_{k-1})/\ac, \Cmax \} - \Gmax$, then
\begin{align*}
S_k & \leq S_{k-1} + \ac (G_{k-1} + D_{k-1} + \D_{k-1}) - \frac{1}{\ad} D_{k-1} \\
& \leq S_{k-1} + \ac (\Gmax + \D_{k-1}) = \Sl_k.
\end{align*}

\item If $- \min\{ \ad S_{k-1}, \Dmax \} - \Gmax \leq \D_{k-1} < -\Gmax$, then
\begin{align*}
S_k & \leq S_{k-1} + \ac C_{k-1} - \frac{1}{\ad} ( -G_{k-1} + C_{k-1} - \D_{k-1} ) \\
& \leq S_{k-1} - \frac{1}{\ad} (- \Gmax - \D_{k-1}) = \Sl_k.
\end{align*}

\item If $\D_{k-1} < -\Gmax - \min\{ \ad S_{k-1}, \Dmax\}$, then $S_k = \Sl_k = S_{k-1} - \min\{S_{k-1}, \Dmax/\ad\}$.
\end{enumerate}

By induction we show that the stationary policy $\pl$ is optimal for the finite horizon stochastic program. Then for any $n \in \N$ and policy $\pi$, we have
\begin{align*}
\frac{1}{n} \Jl_1(\pl, S_1) \leq \frac{1}{n} \Jl_1(\pi, S_1),
\end{align*} 
and thus
\begin{align*}
\Jl(\pl, S_1) = \limsup_{n\to\infty} \frac{1}{n} \Jl_1(\pl,S_1) 
\leq \limsup_{n\to\infty} \frac{1}{n} \Jl_1(\pi,S_1) = \Jl(\pi,S_1).
\end{align*}

\endproof

% residual 
%--------------------
\section{Proofs for the residual power system results}
\label{app:residual}

% lemma for bounding expected stored power
%--------------------
We first establish bounds on the expected stored power.

\begin{lemma}
\label{lm:storage-bounds}
Suppose that $S_{i+1} = \max\{ S_i + \T_i, 0 \}$ for $i = 1,2,\ldots,n$, where $S_1 \geq 0$ is given and $\T_i$, $i = 1,2,\ldots,n$, is a sequence of IID random variables with mean $\mu_\T$ and variance $\s_\T^2$. 
\begin{enumerate}
\item If $\mu_\T \leq 0$, then
$\E[S_{k+1}] \leq \sqrt{k(\mu_\T^2 + \s_\T^2) + S_1^2 }$.

\item If $\mu_\T > 0$, then
$S_1 + k\mu_\T \leq \E[S_{k+1}] \leq k \mu_\T + \sqrt{k}\log_2(4k) \sqrt{\s_\T^2 + S_1^2/k}$.
\end{enumerate}
\end{lemma}

% proof of lemma
%--------------------
\proof{}
For $\mu_\T \leq 0$, consider
\begin{align*}
\E[S_{k+1}^2] & \leq \E[(S_k + \T_k)^2] \\
& = \E[S_k^2] + 2\E[S_k]\E[\T_k] + \E[\T_k^2] \\
& \leq \E[S_k^2] + (\mu_\T^2 + \s_\T^2) \\
& \leq S_1^2 + k(\mu_\T^2 + \s_\T^2),
\end{align*}
where the equality follows by the independence between $\T_k$ and $(\T_1, \T_2, \ldots, \T_{k-1})$ and the last inequality follows by induction. Thus, by Jensen's inequality,
\begin{align*}
\E[S_{k+1}] \leq \sqrt{\E[S_{k+1}^2]} \leq \sqrt{k(\mu_\T^2 + \s_\T^2) + S_1^2 }.
\end{align*}

For $\mu_T > 0$, the expected stored power can be lower bounded as
\begin{align*}
\E[S_{k+1}] \geq \E[S_k + \T_k] \geq S_1 + k\mu_\T,
\end{align*}
where the last inequality follows by induction. Next, we consider
\begin{align*}
\E[ S_{k+1} - k\mu_\T ]
& = \E[ \max\{ S_k + \T_k, 0 \} - k\mu_\T ] \\
& = \E\left[ \max\left\{ 0,\ S_1 + \sum_{i=1}^k \T_i,\ \max_{2 \leq j \leq k} \sum_{i = j}^k \T_i \right\} - k\mu_\T \right] \\
& \leq \E\left[ \max_{1 \leq j \leq k} \left| \sum_{i = j}^k \TT_i \right| \right],
\end{align*}
where $\TT_i = \T_i - \mu_\T$ for $i = 2,3,\ldots,k$ and $\TT_1 = S_1 + \T_1 - \mu_\T$. Note that $\E[\T_i\T_j] = 0$ for $i \neq j$. By the inequality in
%%%(Doob 1953, Lemma 4.1)
\citet[Lemma 4.1]{Doob1953}, 
we have
\begin{align*}
\E\left[ \max_{1 \leq j \leq k} \left| \sum_{i = j}^k \TT_i \right|^2 \right] \leq \log_2(4k)^2 (k\s_\T^2 + S_1^2).
\end{align*}
Therefore,
\begin{align*}
\E[S_{k+1}] & \leq \E\left[ \max_{1 \leq j \leq k} \left| \sum_{i = j}^k \TT_i \right| \right] + k \mu_\T
\leq \sqrt{ \E\left[ \left( \max_{1 \leq j \leq k} \left| \sum_{i = j}^k \TT_i \right| \right)^2 \right] } + k \mu_\T \\
& = \sqrt{ \E\left[ \max_{1 \leq j \leq k} \left| \sum_{i = j}^k \TT_i \right|^2 \right] } + k \mu_\T
\leq \sqrt{k}\log_2(4k)\sqrt{\mu_\T^2 + S_1^2/k} + k \mu_\T,
\end{align*}
where the second inequality follows by Jensen's inequality.
\endproof

% asymptotic fast-ramping generation
%--------------------
\proof{Proof of Propositions~\ref{prop:fastramp-asymp} and~\ref{prop:fastramp-over}.}
Since Proposition~\ref{prop:fastramp-asymp} is a special case of Proposition~\ref{prop:fastramp-over}, we only need to prove the general case.
For unlimited $\Gmax$ and $\Smax$, the optimal policy $\pg$ satisfies
\begin{align*}
0 & = G_i - C_i + D_i + \D_i, \\
S_{i+1} & = S_i + \ac C_i - \frac{1}{\ad} D_i,
\end{align*}
for $i = 1,2,\ldots$. Note that $C_i = \D_i^+$. Now we consider
\begin{align*}
\frac{1}{n} \sum_{i=1}^n G_i & = \frac{1}{n} \sum_{i=1}^n ( C_i - D_i - \D_i ) \\
& = \frac{1}{n} \sum_{i=1}^n \left( (1-\a)C_i - \D_i \right) + \frac{1}{n} \ad ( S_{n+1} - S_1 ) \\
& = \frac{1}{n} \sum_{i=1}^n \left( (1-\a) \D_i^+ - \D_i \right) + \frac{1}{n} \ad ( S_{n+1} - S_1 )\\
& = \frac{1}{n} \sum_{i=1}^n \left( \D_i^- -\a \D_i^+ \right) + \frac{1}{n} \ad (S_{n+1} - S_1).
\end{align*}
Next let $\T_i = \ac \D_i^+ - \D_i^-/\ad$ for $i = 1,2,\ldots$, which is sequence of IID random variables with mean $\mu_\T = \E[\T_i]$ and variance $\s_\T^2 = \Var[\T_i]$. 
Suppose that $\mu_\T \geq 0$. We only need to show that $\limsup_{n \to \infty} \E[S_{n+1}/n] = 0$ since
\begin{align*}
\limsup_{n\to\infty} \E\left[ \frac{1}{n} \sum_{i=1}^n \left( \D_i^- - \a \D_i^+ \right) - \frac{1}{n} \ad S_1 \right] = \E[\D^- - \a\D^+].
\end{align*}
By Lemma~\ref{lm:storage-bounds}, we have
\begin{align*}
0 & \leq \limsup_{n\to\infty} \E\left[ \frac{1}{n} S_{n+1} \right] 
\leq \limsup_{n\to\infty} \frac{1}{n} \sqrt{n(\mu_\T^2+\s_\T^2) + S_1^2} = 0.
\end{align*}
Suppose that $\mu_\T < 0$. By Lemma~\ref{lm:storage-bounds}, we have
\begin{align*}
\limsup_{n\to\infty} \E\left[ \frac{1}{n} \sum_{i=1}^n G_i \right]
& \leq \E\left[ \D^- -\a \D^+ \right] + \limsup_{n\to\infty} \frac{1}{n} \ad \left( n \mu_\T + \sqrt{n}\log_2(4k)\sqrt{\s_\T^2+S_1^2/n} \right) \\
& = \E\left[ \D^- -\a \D^+ \right] + \ad\mu_\T = 0.
\end{align*}
\endproof

% laplace fast-ramping generation
%--------------------
\proof{Proof of Proposition~\ref{prop:fastramp-lapl}.}
To find the expected average fast-ramping generation $\Jcg(\pg,S_1)$ under the Laplace assumption, we need to show that there exist a constant $\eta$ and a function $v(s)$ such that
\begin{align}
\label{equ:acoe}
\eta + v(s) & = \E\left[ \Gg + v\left(s + \ac \Cg - \frac{1}{\ad} \Dg \right)\right],
\end{align}
and $\lim_{n \to \infty} \E\left[ v(S_n)/n \right] = 0$
%%%(Arapostathis et al. 1993)
\cite{Arapostathis1993}. 
Then the average cost is equal to $\eta$. Now we verify that for the policy $\pg$ in Theorem~\ref{thm:policy-generation},
\begin{align*}
\eta & = \frac{1}{2\l} \left( \frac{(1-\a)(1-e^{-\l\Gmax})}{1-\a e^{-(1/\ac-\ad)\l\Smax/2}} \right), \\
v(s) & = - \left( \frac{\ad(1-e^{-\l\Gmax})}{1 - \a e^{-(1/\ac-\ad)\l\Smax/2}} \right) s + \frac{1}{\l} \left( \frac{1+\a}{1-\a} \right) \left( \frac{\a e^{-(1/\ac-\ad)\l(\Smax-s)/2}}{1 - \a e^{-(1/\ac-\ad)\l\Smax/2}} \right) (1-e^{-\l\Gmax})
\end{align*}
satisfy
\begin{align*}
\eta + v(s) & = \E\left[\Gg + v\left(s + \ac \Cg - \frac{1}{\ad} \Dg \right)\right] \\
& = \E\left[ v(\Smax); \frac{1}{\ac}(\Smax -s) \leq \D \right] + \E\left[  v(s+\ac\D); 0 \leq \D <  \frac{1}{\ac}(\Smax-s) \right] \\
&\qquad + \E\left[  v\left(s+ \frac{1}{\ad}\D\right); -\ad s \leq \D < 0 \right] + \E\left[  - \D - \ad s + v(0); -\Gmax - \ad s \leq \D < - \ad s \right] \\
&\qquad + \E\left[ \Gmax + v(0); \D < - \Gmax - \ad s \right],
\end{align*}
where for random variable $X$ and a set $A$ we define $\E[X;A] = \E[X|A]\P(A)$. Furthermore, $v(s)$ is bounded for $s \in [0,\Smax]$, and thus $\lim_{n \to \infty} \E\left[ v(S_n)/n \right] = 0$.
\endproof

% stationary distribution
%--------------------
\proof{Proof of Proposition~\ref{prop:fastramp-dist}.}
Since the optimal policy $\pg$ in Theorem~\ref{thm:policy-generation} is stationary, the corresponding stored power sequence is a Markov process
\begin{align*}
S_{i+1} & = \begin{cases}
\Smax & \text{if}\ (\Smax-S_i)/\ac \leq \D_i \\
S_i + \ac \D_i & \text{if}\ 0 \leq \D_i < (\Smax - S_i)/\ac \\
S_i + \D_i/\ad & \text{if}\ -\ad S_i \leq \D_i < 0 \\
0 & \text{if}\ \D_i < - \ad S_i.
\end{cases}
\end{align*}
Let $F_i$ be the cdf of $S_i$. Then for $0 \leq s < \Smax$,
\begin{align}
\label{equ:fs-recur}
F_{i+1}(s) & = \int_{-\infty}^{-\ad(\Smax-s)} \frac{\l}{2} e^{\l \d} F_i(\Smax) d\d
+ \int_{-\ad(\Smax-s)}^0 \frac{\l}{2} e^{\l \d} F_i(s - \d/\ad) d\d  \notag \\
&\qquad +  \int_0^{(1/\ac)s} \frac{\l}{2} e^{-\l \d} F_i(s-\ac \d) d\d.
\end{align}
For $s < 0$, $F_{i+1}(s) = 0$, and for $s \geq \Smax$, $F_{i+1}(s) = 1$.
Let
\begin{align*}
F_S(s) & = 
\begin{cases}
0 & \text{if}\ s < 0 \\
\displaystyle \frac{1}{1-\a e^{-(1/\ac-\ad)\l\Smax/2}}\left( 1 - \frac{1+\a}{2} e^{-(1/\ac-\ad)\l s/2} \right) & \text{if}\ 0 \leq s < \Smax \\
1 & \text{if}\ s \geq \Smax
\end{cases}.
\end{align*}
It can be verified that~\eqref{equ:fs-recur} is satisfied with $F_i(s) = F_{i+1}(s) = F_S(s)$ for all $s$.
Now we only need to show that the Markov chain is irreducible with respect to $F_S(s)$, which implies that the stationary distribution is unique
%%%(Gilks et al. 1995, Theorem 4.1)
\cite[Theorem~4.1]{Gilks1995}. 
Let $S_i = s$, and let $\Bc = (b_1,b_2) \subset [s,\Smax]$ such that $\P\{S \in \Bc\} > 0$ under the stationary distribution, that is, $b_2 > b_1$. Then 
\begin{align*}
\P\{S_{i+1} \in \Bc\} & = \int_{(1/\ac)(b_1-s)}^{(1/\ac)(b_2-s)} \frac{\l}{2} e^{ - \l\d} d\d
= \frac{1}{2}\left( e^{-(1/\ac)\l(b_1-s)} - e^{-(1/\ac)\l(b_2-s)} \right) > 0.
\end{align*}
If $\Bc = \{0\}$, then
\begin{align*}
P\{S_{i+1} \in \Bc\} = \int_{-\infty}^{-\ad s} \frac{\l}{2} e^{\l\d} d\d 
= \frac{1}{2} e^{-\ad\l s} > 0.
\end{align*}
Similarly, $\P\{S_{i+1} \in \Bc\} > 0$ for $\Bc = (b_1,b_2) \subset [0,s]$ and $\{\Smax\}$. Thus, we can generalize $\Bc$ to any set such that $\P\{S \in \Bc\} > 0$. Therefore, the Markov chain is irreducible.

Using stationary distribution $F_S(s)$ above, we can find the corresponding distribution of the fast-ramping generation
\begin{align*}
F_G(g) & = 1 - \int_{-\infty}^{-g-\ad\Smax} \frac{\l}{2} e^{\l\d} d\d 
- \int_{-g-\ad\Smax}^{-g} \frac{\l}{2}e^{\l\d} F_S\left(-\frac{1}{\ad}(\d+g)\right) d\d \\
& = 1 - \frac{1}{2} e^{-\l g -\l\ad\Smax} \\
&\qquad - \frac{1}{1-\a e^{-(1/\ac-\ad)\l\Smax/2}} \left(  \frac{1-\a}{2}e^{-\l g} - \frac{1}{2}e^{-\l g-\l\ad\Smax} + \frac{\a}{2} e^{-\l g-(1/\ac+\ad)\l\Smax/2}\right) \\
& = 1 - \frac{1-\a}{2(1-\a e^{-(1/\ac-\ad)\l\Smax/2})} e^{-\l g}
\end{align*}
for $g \geq 0$, $F_G(g) = 0$ for $g < 0$, and $F_G(g) = 1$ for $g \geq \Gmax$.
\endproof

% asymptotic loss of load probability
%--------------------
\proof{Proof of Proposition~\ref{prop:lossofload-asymp}.}
For any $u_{k+1} > 0$, consider
\begin{align*}
\E\left[1_{\{\D_{k+1} < -\Gmax - \ad S_{k+1}\}}\right]
& = \E\left[1_{\{\D_{k+1} < -\Gmax - \ad S_{k+1}\}};\ S_{k+1} \leq u_{k+1} \right] \\
&\qquad + \E\left[1_{\{\D_{k+1} < -\Gmax - \ad S_{k+1}\}}; S_{k+1} > u_{k+1} \right] \\
& \leq \P\{ S_{k+1} \leq u_{k+1} \} + \E\left[1_{\{\D_{k+1} < -\Gmax - \ad u_{k+1}\}}\right] \\
& = \P\{ S_{k+1} \leq u_{k+1} \} + F_\D( -\Gmax - \ad u_{k+1} )
\end{align*}
and choose $u_{k+1} = \E[S_{k+1}]/2$. Then we have
\begin{align*}
\E\left[1_{\{\D_{k+1} < -\Gmax - \ad S_{k+1}\}}\right]
& \leq \P\left\{ S_{k+1} < \frac{1}{2} \E[S_{k+1}] \right\} + F_\D\left(-\Gmax - \frac{\ad\E[S_{k+1}]}{2}\right) \\
& = \P\left\{ \E[S_{k+1}] - S_{k+1} > \frac{1}{2} \E[S_{k+1}] \right\} + F_\D\left(-\Gmax - \frac{\ad\E[S_{k+1}]}{2}\right).
\end{align*}
For unlimited $\Smax$, under the optimal policy $\pl$ in Theorem~\ref{thm:policy-lossofload}, we have
\begin{align*}
S_{k+1} & = \max\left\{ S_k + \T_k, 0 \right\}, 
\end{align*}
where $\T_k = \ac(\Gmax+\D_k)^+ - (1/\ad)(\Gmax+\D_k)^-$ with $\mu_\T = \E[\T_k] > 0$ and $\s_\T^2 = \Var[\T_k]$. Then by Lemma~\ref{lm:storage-bounds},
\begin{align*}
S_1 + k\mu_\T \leq \E[S_{k+1}] \leq k \mu_\T + \sqrt{k}\log_2(4k) \sqrt{\s_\T^2 + S_1^2/k}.
\end{align*}
Thus, 
\begin{align*}
\E\left[1_{\{\D_{k+1} < -\Gmax - \ad S_{k+1}\}}\right]
& \leq \P\left\{ \left(\E[S_{k+1}] - S_{k+1}\right)^2 > \frac{1}{4} \E[S_{k+1}]^2 \right\} + F_\D\left(-\Gmax - \frac{\ad}{2}(k\mu_\T+S_1) \right) \\
& \leq \frac{4\Var(S_{k+1})}{\E[S_{k+1}]^2} + F_\D\left(-\Gmax - \frac{\ad}{2}(k\mu_\T+S_1) \right) \\
& \leq \frac{4\Var(S_{k+1})}{(k\mu_\T+S_1)^2} + F_\D\left(-\Gmax - \frac{\ad}{2}(k\mu_\T+S_1) \right),
\end{align*}
where the second inequality follows by the Chebyshev's inequality. Next we consider
\begin{align*}
\Var[S_{k+1}] & = \E[S_{k+1}^2] - \left(\E[S_{k+1}] \right)^2 \\
& \leq \E\left[ (S_k + \T_k)^2 \right] - (k\mu_\T+S_1)^2 \\
& = \E\left[S_k^2 \right] + 2 \mu_\T \E[S_k] + (\mu_\T^2+\s_\T^2) - (k\mu_\T+S_1)^2 \\
& \leq \E\left[S_k^2 \right] + 2 \mu_\T \left( (k-1) \mu_\T + \sqrt{k-1}\log_2(4(k-1)) \sqrt{\s_\T^2 + \frac{S_1^2}{k-1}} \right) + (\mu_\T^2+\s_\T^2) - (k\mu_\T+S_1)^2 \\
& \leq \sum_{j=2}^k 2 \mu_\T \left( (j-1)\mu_\T + \sqrt{j-1}\log_2(4(j-1)) \sqrt{\s_\T^2 + S_1^2} \right) + k(\mu_\T^2+\s_\T^2) + S_1^2 - (k\mu_\T+S_1)^2 \\
& = \sum_{j=2}^k 2 \mu_\T \sqrt{j-1}\log_2(4(j-1)) \sqrt{\s_\T^2 + S_1^2} + k(k-1)\mu_\T^2 + k(\mu_\T^2+\s_\T^2) + S_1^2 - (k\mu_\T+S_1)^2 \\
& = \sum_{j=2}^k 2 \mu_\T \sqrt{j-1}\log_2(4(j-1)) \sqrt{\s_\T^2 + S_1^2} + k (\s_\T^2 - 2\mu_\T S_1),
\end{align*}
where the third inequality follows by induction. Since
\begin{align*}
\sum_{j=2}^k 2 \mu_\T\sqrt{\s_\T^2 + S_1^2} \left( \sqrt{j-1} \frac{\ln 4(j-1) }{\ln 2} \right)
& \leq 2 \mu_\T\sqrt{\s_\T^2 + S_1^2} \int_1^k \sqrt{x} \frac{\ln 4x}{\ln 2} dx \\
& = 2 \mu_\T\sqrt{\s_\T^2 + S_1^2} \left( \left. \frac{2}{3\ln 2} x^{3/2} \ln 4x \right|_1^k - \int_1^k \frac{2}{3\ln 2} x^{1/2} dx \right) \\
& \leq c_1 k^{5/3},
\end{align*}
for some constant $c_1 > 0$ and for sufficiently large $k$, we have
\begin{align*}
\limsup_{n\to\infty}\frac{1}{n} \sum_{k=1}^n \E\left[1_{\{\D_{k+1} < -\Gmax - \ad S_{k+1}\}}\right]
& \leq \limsup_{n\to\infty}\frac{1}{n} \sum_{k=1}^n \left( \frac{4\Var(S_{+1})}{(k\mu_\T+S_1)^2} + F_\D\left(-\Gmax - \frac{\ad}{2}(k\mu_\T+S_1) \right) \right) \\
& \leq \limsup_{n\to\infty}\frac{1}{n} \left( c_2 + c_3 \sum_{k=k'}^n \left( \frac{1}{k^{1/3}} + \frac{1}{k}\right) \right) = 0,
\end{align*}
where $k'$, $c_2$, and $c_3$ are constants.
\endproof

% laplace loss of load probability reducing rate
%--------------------
\proof{Proof of Proposition~\ref{prop:lossofload-rate}.}
We first establish a lower bound on the expected loss of load probability. For power storage capacity $\Smax$, $\Jcl(\pl,S_1)$ is lower bounded by the expected loss of load probability associated with the stored power sequence $S_i=\Smax$ for $i=1,2,\ldots$, that is,
\begin{align*}
\Jcl(\pl,S_1) & \geq \limsup_{n \to\infty} \E\left[ \frac{1}{n} \sum_{i=1}^n 1_{\{\D_i < - \Gmax - \ad\Smax\}}\right] = \frac{1}{2} e^{-\l\Gmax-\ad\l\Smax}.
\end{align*}
Thus,
\begin{align*}
\lim_{\Smax \to \infty} \frac{\ln \Jcl(\pl,S_1)}{\Smax} \geq -\ad\l =  \gamma_{\mathrm{min}}.
\end{align*}

Next, we establish an upper bound on $\Jcl(\pi,S_1)$ by considering a suboptimal policy $\pi$ in Table~\ref{table:lossofload-rate}. Let $F_i$ be the cdf of $S_i$. Then for $0 \leq s < \Smax$,
\begin{align*}
F_{i+1}(s) & = \int_{-\infty}^{-\Gmax-\ad(\Smax-s)} \frac{\l}{2} e^{\l\d} d\d
+ \int_{-\Gmax-\ad(\Smax-s)}^{-\Gmax} \frac{\l}{2} e^{\l\d} F_i\left( s - \frac{\d+\Gmax}{\ad} \right) d\d \\
&\qquad + \int_{-\Gmax}^0 \frac{\l}{2} e^{\l\d} F_i(s) d\d 
+ \int_0^{(1/\ac)s} \frac{\l}{2} e^{-\l\d} F_i(s-\ac\d) d\d \\
& = \frac{1}{2} e^{-\l\Gmax-\ad\l(\Smax-s)} + \frac{1}{2}(1-e^{-\l\Gmax}) F_i(s) \\
&\qquad + \int_{-\Gmax-\ad(\Smax-s)}^{-\Gmax} \frac{\l}{2} e^{\l\d} F_i\left( s - \frac{\d+\Gmax}{\ad} \right) d\d 
+ \int_0^{(1/\ad)s} \frac{\l}{2} e^{-\l\d} F_i(s - \ac\d) d\d.
\end{align*}
It can be verified that
\begin{align*}
F_S(s) & = \frac{e^{-\l\Gmax}}{\a e^{\l_0\Smax} - e^{-\l\Gmax}} \left( -1 + \frac{1+\a}{1+e^{-\l\Gmax}} e^{\l_0 s}\right)
\end{align*}
for $0 \leq s < \Smax$ is a stationary distribution of the process $\{S_i\}$, where
\begin{align*}
0 & < \l_0 = \frac{\a - e^{-\l\Gmax}}{\ac(1+e^{-\l\Gmax})} \l < \ad \l.
\end{align*}
By similar steps in the proof of Proposition~\ref{prop:fastramp-dist}, we can show that the stationary distribution is unique. Thus, 
\begin{align*}
\Jcl(\pl,S_1) & \leq \int_{-\infty}^{-\Gmax-\ad\Smax} \frac{\l}{2} e^{\l\d} d\d + \int_{-\Gmax-\ad\Smax}^{-\Gmax} \frac{\l}{2} e^{\l\d} F_S\left( - \frac{\d+\Gmax}{\ad} \right) d\d \\
& = \frac{1}{2}e^{-\l\Gmax}\left( \frac{\a - e^{-\l\Gmax}}{e^{\l_0\Smax}-e^{-\l\Gmax}} \right).
\end{align*}
Therefore,
\begin{align*}
\lim_{\Smax \to \infty} \frac{\ln \Jcl(\pl,S_1)}{\Smax} \leq -\l_0 = \gamma_{\mathrm{max}}.
\end{align*}

\begin{table}[h]
\TABLE
{Suboptimal policy for minimizing loss of load probability.
\label{table:lossofload-rate}}
{\begin{tabular}{c|ccc|c}
$\pi_i$ & $G_i$ & $C_i$ & $D_i$ & $S_{i+1}$ \\
\hline
$\displaystyle \frac{\Smax-S_i}{\ac} \leq \D_i$ & $0$ & $\displaystyle \frac{\Smax-S_i}{\ac}$ & $0$ & $\Smax$ \\
$\displaystyle 0 \leq \D_i < \frac{\Smax-S_i}{\ac}$ & $0$ & $\D_i$ & $0$ & $S_i + \ac \D_i$\\
$-\Gmax \leq \D_i < 0$ & $-\D_i$ & $0$ & $0$ & $S_i$ \\
$- \Gmax - \ad S_i \leq \D_i < -\Gmax$ & $-\Gmax$ & $0$ & $-\D_i-\Gmax$ & $\displaystyle S_i + \frac{\Gmax+\D_i}{\ad}$ \\
$\D_i < - \Gmax - \ad S_i$ & $\Gmax$ & $0$ & $\ad S_i$ & $0$ \\
\hline
\end{tabular}}
{}
\end{table}

\endproof

% two-threshold
%--------------------
\section{Proofs for the two-threshold policy}
\label{app:threshold}

% optimality for 2-slot dynamic program
%--------------------
\proof{Proof of Proposition~\ref{prop:threshold}.}
Define the cost-to-go function of the two-slot dynamic program
\begin{align*}
J_k(\pi,S_k) = \E\left[ \left. \sum_{i=k}^2 \rho_1 G_i + \rho_2 1_{\{\D_i < -\Gmax - \ad S_i\}} \right| S_k \right] 
\end{align*}
for $k = 1,2$. When $k=2$, the expected loss of load probability is equal to
\begin{align*}
\E\left[ \left. \rho_2 1_{\{\D_2 < -\Gmax - \ad S_2\}} \right| S_2 \right]
\end{align*}
and does not depend on the control $\pi_2 = (G_2,C_2,D_2)$. Thus, the optimal policy at time $2$ is $\pi_2^* = \pg_2$, that is, ${\Sc}_2 = {\Sd}_2 = 0$. Then the minimum cost-to-go function is
\begin{align*}
J_2(\pi^*,S_2) 
& = \E\left[ \rho_1 \Gg_2 + \rho_2 1_{\{\D_2 < -\Gmax - \ad S_2\}} \right] \\
& = \int_{-\infty}^{-\Gmax - \ad S_2} (\rho_1 \Gmax + \rho_2) f_\D(\d_2)\ d\d_2 - \int_{-\Gmax-\ad S_2}^{-\ad S_2} \rho_1 \d_2 f_\D(\d_2)\ d\d_2.
\end{align*}
Note that $J_2(\pi^*,S_2)$ is convex and decreasing in $S_2$ since
\begin{align*}
\frac{\partial J_2(\pi^*,S_2)}{\partial S_2} 
& = -\ad (\rho_1 \Gmax + \rho_2) f_\D(-\Gmax - \ad S_2) \\
&\qquad + \ad \rho_1 \Gmax f_\D(-\Gmax - \ad S_2) - \int_{-\Gmax-\ad S_2}^{-\ad S_2} \ad \rho_1 f_\D(\d_2)\ d\d_2 \\
& = -\ad \rho_2 f_\D(-\Gmax - \ad S_2) - \int_{-\Gmax}^{0} \ad \rho_1 f_\D(\d_2-\ad S_2)\ d\d_2 \leq 0,
\end{align*}
and $f_\D$ in increasing on $(-\infty,0]$. 

Next we consider the cost-to-go function at time $1$. For a policy $\pi_1$, we consider the following three cases.
\begin{enumerate}
\item If $\D_1 \geq (\Smax-S_1)/\ac$, then
\begin{align*}
\rho_1 G_1 + \rho_2 1_{\{\D_1 < - \Gmax - \ad S_1 \}} + J_2\left(\pi^*,S_1 + \ac C_1 - \frac{1}{\ad} D_1 \right)
& \geq J_2\left(\pi^*, \Smax \right),
\end{align*}
where the lower bound is achieved by $G_1 = 0$, $C_1 = (\Smax-S_1)/\ac$, and $D_1 = 0$.

\item If $\D_1 < -\Gmax -\ad S_1$, then loss of load occurs, and $G_1 = \Gmax$, $C_1 = 0$, and $D_1 = \ad S_1$.

\item For $-\Gmax - \ad S_1 \leq \D_1 < (\Smax - S_1)/\ac$, consider the lower bounds on fast-ramping generation
\begin{align*}
G_1 & \geq \frac{1}{\ac}\left(S_2 - S_1 + \frac{1}{\ad} D_1 \right) - D_1 - \D_1 
\geq - \D_1 + \frac{1}{\ac}(S_2 - S_1), \\
G_1 & \geq C_1 - \ad(\ac C_1 + S_1 - S_2) - \D_1 
\geq - \D_1 + \ad(S_2-S_1).
\end{align*}
Thus, we have
\begin{align*}
\rho_1 G_1 + J_2(\pi^*,S_2) & \geq \max\{J_{1\mathrm{c}}(S_1,S_2),J_{1\mathrm{d}}(S_1,S_2)\},
\end{align*}
where
\begin{align*}
J_{1\mathrm{c}}(S_1,S_2) & = \rho_1 \left(-\d_1 + \frac{1}{\ac}(S_2-S_1) \right) + J_2(\pi^*,S_2), \\
J_{1\mathrm{d}}(S_1,S_2) & = \rho_1 \left(-\d_1 + \ad(S_2-S_1) \right) + J_2(\pi^*,S_2).
\end{align*}
Note that the loss of load probability $\rho_2 1_{\{\D_1 < -\Gmax - \ad S_1 \} }$ does not depend on the control. 
By the convexity of $J_2(\pi^*,S_2)$, $J_{1\mathrm{c}}(S_1,S_2)$ and $J_{1\mathrm{d}}(S_1,S_2)$ are convex. Let 
\begin{align*}
{\Sc}_1 & = \argmin_{0 \leq S_2 \leq \Smax} J_{1\mathrm{c}}(S_1,S_2) = \left( \sup\left\{ 0 \leq S_2 \leq \Smax: \frac{\partial J_2(\pi^*,S_2)}{\partial S_2} \leq - \frac{\rho_1}{\ac}  \right\}  \right)^+, \\
{\Sd}_1 & = \argmin_{0 \leq S_2 \leq \Smax} J_{1\mathrm{d}}(S_1,S_2) = \left( \sup\left\{ 0 \leq S_2 \leq \Smax: \frac{\partial J_2(\pi^*,S_2)}{\partial S_2} \leq - \ad \rho_1 \right\}  \right)^+,
\end{align*}
where $\sup \emptyset = -\infty$.
Then $0 \leq {\Sc}_1 \leq {\Sd}_1 \leq \Smax$. If $\D_1 \geq -\Gmax$, then $0 \leq S_2 \leq \min\{S_1+\ac(\Gmax+\d), \Smax\}$, and if $\D_1 < -\Gmax$, then $0 \leq S_2 \leq \max\{S_1 + (\Gmax+\d)/\ad,0\}$. Therefore, the two-threshold policy with parameters $({\Sc}_1,{\Sd}_1)$ is optimal at time $1$ by the convexity of $J_{1\mathrm{c}}$ and $J_{1\mathrm{d}}$.
\end{enumerate}
\endproof

\end{APPENDICES}

%%
%\theendnotes

% Acknowledgments here
\ACKNOWLEDGMENT{The authors would like to thank Javad Lavaei for feedback that helped simplify and generalize the proofs of Theorems~\ref{thm:policy-generation} and ~\ref{thm:policy-lossofload}. We also thank Ram Rajagopal and Benjamin Van Roy for valuable comments and suggestions that greatly improved the presentation of this paper.}

% References here (outcomment the appropriate case)

% CASE 1: BiBTeX used to constantly update the references
%   (while the paper is being written).
\bibliographystyle{ormsv080} % outcomment this and next line in Case 1
\bibliography{storage-or} % if more than one, comma separated

% CASE 2: BiBTeX used to generate mypaper.bbl (to be further fine tuned)
%\input{mypaper.bbl} % outcomment this line in Case 2

%If you don't use BiBTex, you can manually itemize references as shown below.

%%%%%%%%%%%%%%%%%
\end{document}